\documentclass[aop,MSNbibl,seceqn,dvips]{arximspdf}
\usepackage{accents}
\usepackage{graphicx}
\doi{10.1214/12-AOP797} 
\volume{41}
\issue{5}
\pubyear{2013}
\firstpage{3201}
\lastpage{3240}

\makeatletter
\newcommand{\rrvert}{\vert}
\newcommand{\llvert}{\vert}
\renewcommand{\mathring}[1]{\accentset{\circ}{#1}}
\def\cal{\mathcal}

\newtheorem{Theorem}{Theorem}[section]
\newtheorem{Proposition}[Theorem]{Proposition}
\newproclaim{Assumption}[Theorem]{Assumption}
\newtheorem{Algorithm}{Algorithm}
\newtheorem{Lemma}[Theorem]{Lemma}
\newproclaim{Remark}[Theorem]{Remark}
\newproclaim{nota}{Notation}
\def\E{\mathbb{E}}
\def\F{\mathbb{F}}
\def\N{\mathbb{N}}
\def\P{\mathbb{P}}
\def\Q{\mathbb{Q}}
\def\R{\mathbb{R}}
\def\Z{\mathbb{Z}}

\def\a{\alpha}
\def\ub{ \bar{u} }
\def\O{\Omega}
\def\Ob{\overline{\O}}
\def\om{\omega}
\def\eps{\varepsilon}
\def\xb{{\mathbf{x}}}
\def\abb{ \bar{a}}
\def\bbb{ \bar{b}}

\def\Bc{\cal{B}}
\def\Ec{\cal{E}}
\def\Fc{\cal{F}}
\def\Lc{\cal{L}}
\def\Mc{\cal{M}}
\def\Nc{\cal{N}}
\def\Pc{\cal{P}}
\def\Tc{\cal{T}}
\def\Vc{\cal{V}}

\def\Fb{\overline{\F}}
\def\Fcb{\overline{\Fc}}
\def\Pcb{\overline{\Pc}}

\def\Ab{\overline{A}}
\def\Jb{\overline{J}}

\def\Pb{\overline{\P}}
\def\Qb{{\bar\Q}}

\def\x{\times}

\def\eps{\varepsilon}

\def\1{\mathbf{1}}

\def\lin{\operatorname{lin}}
\newcommand{\eqref}[1]{(\ref{#1})}
\makeatother

\begin{document}
\begin{frontmatter}

\title{Optimal transportation under controlled stochastic dynamics\thanksref{T1}}
\runtitle{Semimartingale transportation problem}
\thankstext{T1}{Supported by the Chair \textit{Financial Risks} of the
\textit{Risk Foundation} sponsored by Soci\'et\'e
G\'en\'erale, the Chair \textit{Derivatives of the Future} sponsored by
the {F\'ed\'eration Bancaire Fran\c{c}aise}, and
the Chair \textit{Finance and Sustainable Development} sponsored by EDF
and CA-CIB.}

\begin{aug}
\author[A]{\fnms{Xiaolu} \snm{Tan}\corref{}\ead[label=e1]{xiaolu.tan@polytechnique.edu}}
\and
\author[A]{\fnms{Nizar} \snm{Touzi}\ead[label=e2]{nizar.touzi@polytechnique.edu}}
\runauthor{X. Tan and N. Touzi}
\affiliation{Ecole Polytechnique, Paris}
\address[A]{CMAP\\
Ecole Polytechnique, Paris\\
91128, Palaiseau\\
France\\
\printead{e1}\\
\phantom{E-mail:\ }\printead*{e2}} 
\end{aug}

\received{\smonth{7} \syear{2011}}
\revised{\smonth{6} \syear{2012}}

%
\begin{abstract}
We consider an extension of the Monge--Kantorovitch optimal
transportation problem.
The mass is transported along a continuous semimartingale, and the cost
of transportation
depends on the drift and the diffusion coefficients of the continuous
semimartingale.
The optimal transportation problem minimizes the cost among all
continuous semimartingales with given initial and terminal
distributions. Our first main result is an extension of the
Kantorovitch duality to this context. We also suggest a
finite-difference scheme combined with the gradient projection
algorithm to approximate the dual value. We prove the convergence of
the scheme, and we derive a rate of convergence.

We finally provide an application in the context of financial
mathematics, which originally motivated our extension of the
Monge--Kantorovitch problem. Namely, we implement our scheme to
approximate no-arbitrage bounds on the prices of exotic options given
the implied volatility curve of some maturity.
\end{abstract}

%


\end{frontmatter}
%
\section{Introduction}
\label{secintroduction}

In the classical mass transportation problem of Monge--Kantorovich, we
fix at first an initial probability distribution $\mu_0$ and a terminal
distribution $\mu_1$ on $\R^d$. An admissible transportation plan is
defined as a random vector $(X_0, X_1)$ (or, equivalently, a joint
distribution on $\R^d \x\R^d$) such that the marginal distributions
are, respectively, $\mu_0$ and $\mu_1$. By transporting the mass from the
position $X_0(\om)$ to the position $X_1(\om)$, an admissible plan
transports a mass from the distribution $\mu_0$ to the distribution
$\mu_1$. The transportation cost is a function of the initial and final
positions, given by $\E[c(X_0, X_1)]$ for some function $c \dvtx\R^d
\x
\R^d \to\R^+$. The Monge--Kantorovich problem consists in minimizing the
cost among all admissible transportation plans. Under mild conditions,
a duality result is established by Kantorovich, converting the problem
into an optimization problem under linear constraints. We refer to
Villani \cite{Villani} and Rachev and Ruschendorf \cite
{RachevRuschendorf} for this classical duality and the richest
development on the classical mass transportation problem.

As an extension of the Monge--Kantorovitch problem, Mikami and Thieul\-len~\cite{Mikami2006}
introduced the following stochastic mass
transportation mechanism. Let $X$ be an $\R^d$-continuous
semimartingale with decomposition
%
\begin{equation}
\label{eqXMT} X_t = X_0 + \int_0^t
\beta_s \,ds + W_t,
\end{equation}
where $W_t$ is a $d$-dimensional standard Brownian motion under the
filtration $\F^{X}$ generated by $X$. The optimal mass transportation
problem consists in minimizing the cost of transportation defined by
some cost functional $\ell$ along all transportation plans with initial
distribution $\mu_0$ and final distribution $\mu_1$:
\[
V(\mu_0, \mu_1):= \inf\E\int_0^1
\ell(s,X_s,\beta_s) \,ds,
\]
where the infimum is taken over all semimartingales given by \eqref
{eqXMT} satisfying $\P\circ X_0^{-1} = \mu_0$ and $\P\circ X_1^{-1}
= \mu_1$. Mikami and Thieullen \cite{Mikami2006} proved a strong
duality result, thus extending the classical Kantorovitch duality to
this context.

Motivated by a problem in financial mathematics, our main objective is
to extend \cite{Mikami2006} to a larger class of transportation plans
defined by continuous semimartingales with absolutely continuous
characteristics:
\[
X_t = X_0 + \int_0^t
\beta_s \,ds + \int_0^t
\sigma_s \,dW_s,
\]
where the pair process $(\alpha:= \sigma\sigma^T,\beta)$ takes values
in some closed convex subset $U$ of $\R^{d\times d}\times\R^d$, and the
transportation cost involves the drift and diffusion coefficients as
well as the trajectory of $X$.

The simplest motivating problem in financial mathematics is the
following. Let $X$ be the price process of some tradable security, and
consider some path-dependent derivative security $\xi(X_t,t\le1)$.
Then, by the no-arbitrage theory, any martingale measure $\P$ (i.e.,
probability measure under which $X$ is a martingale) induces an
admissible no-arbitrage price $\E^\P[\xi]$ for the derivative security
$\xi$. Suppose further that the prices of all $1$-maturity European
call options with all possible strikes are available. This is a
standard assumption made by practitioners on liquid options markets.
Then, the collection of admissible martingale measures is reduced to
those which are consistent with this information, that is, $c_1(y):=\E^\P[(X_1-y)^+]$
is given for all\vadjust{\goodbreak} $y \in\R$ or, equivalently, the
marginal distribution of $X_1$ under $\P$ is given by $\mu_1[y,\infty
):=-\partial^-c_1(y)$, where $\partial^-c_1$ denotes the left-hand side
derivative of the convex function $c_1$. Hence, a natural formulation
of the no-arbitrage lower and upper bounds is $\inf\E^\P[\xi]$ and
$\sup
\E^\P[\xi]$ with optimization over the set of all probability measures
$\P$ satisfying $\P\circ(X_0)^{-1}=\delta_{x}$ and $\P\circ
(X_1)^{-1}=\mu_1$, for some initial value of the underlying asset price
$X_0=x$. We refer to Galichon, Henry-Labord\`ere and Touzi \cite{ght}
for the connection to the so-called model-free superhedging problem. In
Section~\ref{subsecnumericalexample} we shall provide some
applications of our results in the context of variance options $\xi
=\langle\log X\rangle_1$ and the corresponding weighted variance
options extension.

This problem is also intimately connected to the so-called Skorokhod
Embedding Problem (SEP) that we now recall; see Obloj \cite{Obloj2004}
for a review. Given a one-dimensional Brownian motion $W$ and a
centered $|x|$-integrable probability distribution $\mu_1$ on $\R$, the
SEP consists in searching for a stopping time $\tau$ such that
$W_{\tau
} \sim\mu_1$ and $(W_{t \land\tau})_{t \ge0}$ is uniformly
integrable. This problem is well known to have infinitely many
solutions. However, some solutions have been proved to satisfy some
optimality with respect to some criterion (Az\'ema and Yor \cite
{AzemaYor}, Root~\cite{Root1969} and Rost~\cite{Rost1976}). In order
to recast the SEP in our context, we specify the set $U$, where the
characteristics take values, to $U=\R\times\{0\}$, that is,
transportation along a local martingale. Indeed, given a solution $\tau
$ of the SEP, the process $X_t:= W_{\tau\land{t}/{(1-t)}}$ defines
a continuous local martingale satisfying $X_1 \sim\mu_1$. Conversely,
every continuous local martingale can be represented as time-changed
Brownian motion by the Dubins--Schwarz theorem (see, e.g., Theorem 4.6,
Chapter 3 of Karatzas and Shreve \cite{Karatzas1991}).

We note that the seminal paper by Hobson \cite{Hobson1998} is
crucially based on the connection between the SEP and the above problem
of no-arbitrage bounds for a specific restricted class of derivatives
prices (e.g., variance options, lookback option, etc.).
We refer to Hobson \cite{Hobson2009} for an overview on some specific
applications of the SEP in the context of finance.

Our first main result is to establish the Kantorovitch strong duality
for our semimartingale optimal transportation problem. The dual value
function consists in the minimization of $\mu_0(\lambda_0) - \mu_1(\lambda_1)$
over all continuous and bounded functions $\lambda_1$,
where $\lambda_0$ is the initial value of a standard stochastic control
problem with final cost $\lambda_1$. In the Markovian case, the
function $\lambda_0$ can be characterized as the unique viscosity
solution of the corresponding dynamics programming equation with
terminal condition $\lambda_1$.

Our second main contribution is to exploit the dual formulation for the
purpose of numerical approximation of the optimal cost of
transportation. To the best of our knowledge, the first attempt for the
numerical approximation of an intimately related problem, in the
context of financial mathematics, was initiated by Bonnans and Tan
\cite
{BonnansTan}. In this paper, we follow their approach in the context of
a bounded set of admissible semimartingale characteristics. Our
numerical scheme combines the finite difference scheme and the gradient\vadjust{\goodbreak}
projection algorithm. We prove convergence of the scheme, and we derive
a rate of convergence.
We also implement our numerical scheme and give some numerical experiments.

The paper is organized as follows.
Section~\ref{secprobmonge} introduces the optimal mass transportation
problem under controlled stochastic dynamics.
In Section~\ref{secduality} we extend the Kantorovitch duality to our
context by using the classical convex duality approach. The convex
conjugate of the primal problem turns out to be the value function of a
classical stochastic control problem with final condition given by the
Lagrange multiplier lying in the space of bounded continuous functions.
Then the dual formulation consists in maximizing this value over the
class of all Lagrange multipliers. We also show, under some conditions,
that the Lagrange multipliers can be restricted to the subclass of
$C^\infty$-functions with bounded derivatives of any order. In the
Markovian case, we characterize the convex dual as the viscosity
solution of a dynamic programming equation in the Markovian case in
Section~\ref{secPDEdynprog}. Further, when the characteristics are
restricted to a bounded set, we use the probabilistic arguments to
restrict the computation of the optimal control problem to a bounded
domain of $\R^d$.

Section~\ref{secnumericalapprox} introduces a numerical scheme to
approximate the dual formulation in the Markovian case. We first use
the finite difference scheme to solve the control problem. The
maximization is then approximated by means of the gradient projection
algorithm. We provide some general convergence results together with
some control of the error.
Finally, we implement our algorithm and provide some numerical examples
in the context of its applications in financial mathematics. Namely, we
consider the problem of robust hedging weighted variance swap
derivatives given the prices of European options of all strikes. The
solution of the last problem can be computed explicitly and allows to
test the accuracy of our algorithm.

\begin{nota*}
Given a Polish space $E$, we denote by
$\mathbf{M}(E)$ the space of all Borel probability measures on $E$,
equipped with the weak topology, which is also a Polish space. In
particular, $\mathbf{M}(\R^d)$ is the space of all probability measures
on $(\R^d, \Bc(\R^d))$. $S_d$ denotes the set of $d \x d$ positive
symmetric matrices. Given $u = (a,b) \in S_d \x\R^d$, we define $|u|$
by its $L^2$-norm as an element in $\R^{d^2+d}$. Finally, for every
constant $C \in\R$, we make the convention $\infty+ C = \infty$.
\end{nota*}

\section{The semimartingale transportation problem}
\label{secprobmonge}

Let $\O:= C([0,1], \R^d)$ be the canonical space, $X$ be the
canonical process
\[
X_t(\om):= \om_t \qquad\mbox{for all } t\in[0,1],
\]
and $\F= ( \Fc_t )_{1 \le t \le1}$ be the canonical filtration
generated by $X$. We recall that $\Fc_t$ coincides with the Borel
$\sigma$-field on $\Omega$ induced by the seminorm $| \om|_{\infty,t}:= \sup_{0 \le s \le t} |\om_s|$,
$\om\in\O$ (see, e.g., the
discussions in Section 1.3, Chapter 1 of Stroock and Varadhan \cite
{Stroock1979}).\vadjust{\goodbreak}

Let $\P$ be a probability measure on $(\O,\Fc_1)$ under which the
canonical process~$X$ is a $\F$-continuous semimartingale. Then, we
have the unique continuous decomposition w.r.t. $\F$:
%
\begin{equation}
\label{eqXdecomposition} X_t = X_0 + B^{\P}_t
+ M^{\P}_t,\qquad t\in[0,1], \P\mbox{-a.s.},
\end{equation}
where $B^{\P} = (B^{\P}_t)_{0 \le t \le1}$ is the finite variation
part and $M^{\P} = (M^{\P}_t)_{0 \le t \le1}$ is the local martingale
part satisfying $B_0 = M_0 = 0$. Denote by $A^{\P}_t:= \langle M^{\P}
\rangle_t$ the quadratic variation of $M^{\P}$ between $0$ and $t$ and
$A^{\P} = (A^{\P}_t)_{0 \le t \le1}$. Then, following Jacod and
Shiryaev \cite{Jacod1987}, we say that the $\P$-continuous
semimartingale $X$ has characteristics $(A^{\P},B^{\P})$.

In this paper, we further restrict to the case where the processes
$A^{\P}$ and $ B^{\P}$ are absolutely continuous in $t$ w.r.t. the
Lebesgue measure, $\P$-a.s. Then there are $\F$-progressive processes $
\nu^{\P} = (\alpha^{\P}, \beta^{\P})$ (see, e.g., Proposition I.3.13
of \cite{Jacod1987}) such that
%
\begin{eqnarray}
\label{eqRepresentationAB} A^{\P}_t = \int
_0^t \a^{\P}_s \,ds,\qquad
B^{\P}_t = \int_0^t
\beta^{\P}_s \,ds, \qquad\P\mbox{-a.s.}\mbox{ for all } t
\in[0,1].
\end{eqnarray}

\begin{Remark}\label{remsemimart2Ito}
By Doob's martingale representation theorem (see, e.g., Theorem 4.2 in
Chapter 3 of Karatzas and Shreve \cite{Karatzas1991}), we can find a
Brownian motion~$W^{\P}$ (possibly in an enlarged space) such that $X$
has an It\^o representation:
\[
X_t = X_0 + \int_0^t
\beta^{\P}_s \,ds + \int_0^t
\sigma^{\P}_s \,dW^{\P}_s,
\]
where $\sigma^{\P}_t = (\a^{\P}_t)^{1/2}$ [i.e., $\a^{\P}_t =
\sigma^{\P
}_t (\sigma^{\P}_t)^T$].
\end{Remark}

\begin{Remark}\label{remuniquenu}
With the unique processes $(A^{\P}, B^{\P})$, the progressively
measurable processes $\nu^{\P} = (\a^{\P}, \beta^{\P})$ may not be
unique. However, they are unique in sense $d \P\x \,dt$-a.e. Since the
transportation cost defined below is a $d \P\x \,dt$ integral, then the
choice of $\nu^{\P} = (\a^{\P}, \beta^{\P})$ will not change the cost
value and then is not essential.
\end{Remark}

We next introduce the set $U$ defining some restrictions on the
admissible characteristics:
%
\begin{equation}
\label{eqconvexitysetSigma} U \mbox{ closed and convex subset of }
S_d \x\R^d,
\end{equation}
and we denote by $\Pc$ the set of probability measures $\P$ on $\O$
under which $X$ has the decomposition \eqref{eqXdecomposition} and
satisfies \eqref{eqRepresentationAB} with characteristics $\nu^{\P
}_t:= (\alpha^{\P}_t, \beta^{\P}_t) \in U$, $d\P\x \,dt$-a.e.

Given two arbitrary probability measures $\mu_0$ and $\mu_1$ in $\mathbf{
M}(\R^d)$, we also denote
%
\begin{eqnarray}
\Pc(\mu_0) &:=& \bigl\{ \P\in\Pc\dvtx \P\circ X_0^{-1}
= \mu_0 \bigr\}, \label{Pmu0}
\\
\Pc(\mu_0, \mu_1) &:=& \bigl\{ \P\in\Pc(
\mu_0) \dvtx\P\circ X_1^{-1} = \mu_1
\bigr\}. \label{Pmu0mu1}
\end{eqnarray}

\begin{Remark}
(i) In general, $\Pc(\mu_0, \mu_1)$ may be empty. However, in the
one-dimensional case $d=1$ and $U = \R^+ \x\R$, the initial
distribution $\mu_0 = \delta_{x_0}$ for some constant $x_0 \in\R$, and
the final distribution satisfies $\int_{\R} |x| \mu_1(dx) < \infty$, we
now verify that $\Pc(\mu_0, \mu_1)$ is not empty. First, we can choose
any constant in $\R$ for the drift part, hence, we can suppose, without
loss of generality, that $x_0 = 0$ and $\mu_1$ is centered distributed,
that is, $\int_{\R} x \mu_1(dx) = 0$. Then, given a Brownian motion
$W$, by Skorokhod embedding (see, e.g., Section 3 of Obloj \cite
{Obloj2004}), there is a stopping time $\tau$ such that $W_{\tau}
\sim
\mu_1$ and $(W_{ t \land\tau})_{t \ge0}$ is uniformly integrable.
Therefore, $M = (M_t)_{0 \le t \le1}$ defined by $M_t:= W_{\tau\land
{t}/{(1-t)}} $ is a continuous martingale with marginal distribution
$\P\circ M_1^{-1} = \mu_1$. Moreover, its quadratic variation
$\langle
M \rangle_t = \tau\land\frac{t}{1-t}$ is absolutely continuous in $t$
w.r.t. Lebesgue for every fixed $\om$, which can induce a probability on
$\O$ belonging to $\Pc(\mu_0, \mu_1)$.\vspace*{-6pt}
\begin{longlist}[(ii)]
\item[(ii)] Let $d=1$, $U=\R^+\times\{0\}$, $\mu_0=\delta_{x_0}$ for some
constant $x_0\in\R$, and $\mu_1$ as in~(i) with $\int x\mu_1(dx)=x_0$.
Then, by the above discussion, we also have $\Pc(\mu_0,\mu_1)\neq
\varnothing$.
\end{longlist}
\end{Remark}

The semimartingale $X$ under $\P$ can be viewed as a vehicle of mass
transportation, from the $\P$-distribution of $X_0$ to the $\P
$-distribution of $X_1$. We then associate~$\P$ with a transportation cost
%
\begin{equation}
\label{eqdJalpha} J(\P):=\E^{\P} \int_0^1
L\bigl(t,X, \nu^{\P}_t\bigr) \,dt,
\end{equation}
where $\E^{\P}$ denotes the expectation under the probability measure
$\P$, and $L\dvtx[0,1]\times\O,\times U\longrightarrow\R$. The above
expectation is well defined on $\R^+ \cup\{+\infty\}$ in view of the
subsequent Assumption~\ref{assumLconvex} which states, in particular,
that $L$ is nonnegative.

Our main interest is on the following optimal mass transportation
problem, given two probability measures $\mu_0$, $\mu_1 \in\mathbf
{M}(\R^d)$:
%
\begin{equation}
\label{eqVprimal} V(\mu_0, \mu_1):=
\inf_{\P\in\Pc(\mu_0, \mu_1)} J(\P),
\end{equation}
with the convention $\inf\varnothing= \infty$.

\section{The duality theorem}
\label{secduality}

The main objective of this section is to prove a duality result for
problem \eqref{eqVprimal} which extends the classical Kantorovitch
duality in optimal transportation theory.

This will be achieved by classical convex duality techniques which
require to verify that the function $V$ is convex and lower
semicontinuous. For the general theory on duality analysis in Banach
spaces, we refer to Bonnans and Shapiro \cite{Bonnans2000} and Ekeland
and Temam \cite{Ekeland1999}. In our context, the value function of
the optimal transportation problem is defined on the Polish space of
measures on $\R^d$, and our main reference is Deuschel and Stroock
\cite{Deuschel1989}.

\subsection{The main duality result}

We first formulate the assumptions needed for our duality result.

\begin{Assumption}\label{assumLconvex}
The function $L\dvtx(t, \xb, u)\in[0,1] \x\O\x U \mapsto L(t, \xb, u)
\in\R^+$ is nonnegative, continuous in $(t, \xb, u)$, and convex in $u$.
\end{Assumption}

Notice that we do not impose any progressive measurability for the
dependence of $L$ on the trajectory $\xb$. However, by immediate
conditioning, we may reduce the problem so that such a progressive
measurability is satisfied.

The next condition controls the dependence of the cost functional on
the time variable.

\begin{Assumption}\label{assumLunifcontinuity}
The function $L$ is uniformly continuous in $t$ in the sense that
\[
\Delta_t L(\eps):= \sup\frac{ |L(s,\xb,u) - L(t,\xb,u)| } {
1+L(t,\xb,u)} \longrightarrow0
\qquad\mbox{as } \eps\to0,
\]
where the supremum is taken over all $0 \le s, t \le1$ such that
$|t-s| \le\eps$ and all $\xb\in\O$, $u \in U$.
\end{Assumption}

We finally need the following coercivity condition on the cost functional.

\begin{Assumption}\label{assumLcroissantp}
There are constants $p > 1$ and $C_0 >0$ such that
\begin{eqnarray*}
|u|^p \le C_0\bigl(1 + L(t,\xb,u)\bigr) < \infty \qquad
\mbox{for every } (t,\xb,u) \in[0,1] \x\O\x U.
\end{eqnarray*}
\end{Assumption}

\begin{Remark}
In the particular case $U=\{I_d\} \x\R^d$, the last condition
coincides with Assumption A.1 of Mikami and Thieullen \cite
{Mikami2006}. Moreover, whenever $U$ is bounded, Assumption \ref
{assumLcroissantp} is a direct consequence of Assumption \ref
{assumLconvex}.
\end{Remark}

Let $C_b(\R^d)$ denote the set of all bounded continuous functions on
$\R^d$ and
\begin{eqnarray*}
\mu(\phi):= \int_{\R^d} \phi(x) \mu(dx) \qquad\mbox{for all } \mu
\in\mathbf{ M}\bigl(\R^d\bigr) \mbox{ and } \phi\in\Lc^1(\mu).
\end{eqnarray*}
We define the dual formulation of \eqref{eqVprimal} by
%
\begin{equation}
\label{eqVdual} \Vc(\mu_0, \mu_1):=
\sup_{\lambda_1 \in C_b(\R^d)} \bigl\{ \mu_0(\lambda_0) -
\mu_1(\lambda_1) \bigr\},
\end{equation}
where
%
\begin{eqnarray}
\label{eqlambda0L} \lambda_0 (x):=\inf_{\P\in\Pc(\delta_{x})}
\E^{\P} \biggl[ \int_0^1 L\bigl(s,X,
\nu^{\P}_s\bigr)\,ds + \lambda_1(X_1)
\biggr],
\end{eqnarray}
with $\Pc(\delta_x)$ defined in \eqref{Pmu0}. We notice that $\mu_0(\lambda_0)$ is well defined since $\lambda_0$ takes value in $\R
\cup\{ \infty\}$, is bounded from below and is measurable by the
following lemma.

\begin{Lemma} \label{lemmdynamprogmu0}
Let Assumptions~\ref{assumLconvex} and~\ref{assumLunifcontinuity}
hold true. Then, $\lambda_0$ is measurable w.r.t. the Borel $\sigma
$-field on $\R^d$ completed by $\mu_0$, and
\[
\mu_0(\lambda_0) = \inf_{\P\in\Pc(\mu_0)}
\E^{\P} \biggl[ \int_0^1 L\bigl(s,X,
\nu^{\P}_s\bigr)\,ds + \lambda_1(X_1)
\biggr].
\]
\end{Lemma}

The proof of Lemma~\ref{lemmdynamprogmu0} is based on a measurable
selection argument and is reported at the end of Section~\ref{sect-proofDPP}.
We now state the main duality result.

\begin{Theorem}\label{theodual}
Let Assumptions~\ref{assumLconvex},~\ref{assumLunifcontinuity} and
\ref{assumLcroissantp} hold. Then
\begin{eqnarray*}
V(\mu_0, \mu_1) = \Vc(\mu_0,
\mu_1)\qquad \mbox{for all } \mu_0, \mu_1 \in
\mathbf{M}\bigl(\R^d\bigr),
\end{eqnarray*}
and the infimum is achieved by some $\P\in\Pc(\mu_0,\mu_1)$ for the
problem $V(\mu_0,\mu_1)$ of~\eqref{eqVprimal}.
\end{Theorem}

The proof of this result is reported in the subsequent subsections.

We finally state a duality result in the space $C_b^{\infty}(\R^d)$ of
all functions with bounded derivatives of any order:
%
\begin{eqnarray}
\overline{\Vc}(\mu_0, \mu_1):=\sup_{\lambda_1 \in C_b^{\infty} (\R^d)}
\bigl\{ \mu_0(\lambda_0) - \mu_1(
\lambda_1) \bigr\}.
\end{eqnarray}

\begin{Assumption}\label{assumLunifcontinuityinX}
The function $L$ is uniformly continuous in $\xb$ in the sense that
\[
\Delta_x L(\eps):= \sup\frac{ |L(t,\xb^1,u) - L(t,\xb^2,u)| }{1+L(t,\xb^2,u)} \longrightarrow0,
\qquad\mbox{as } \eps\to0,
\]
where the supremum is taken over all $0 \le t \le1$, $u \in U$ and all
$\xb^1, \xb^2 \in\O$ such that $|\xb^1 - \xb^2|_{\infty} \le
\eps$.
\end{Assumption}

\begin{Theorem}\label{theoweakdual}
Under the conditions of Theorem~\ref{theodual} together with
Assumption~\ref{assumLunifcontinuityinX}, we have $\Vc=\overline
{\Vc}$ on $\mathbf{M}(\R^d)\times\mathbf{M}(\R^d)$.
\end{Theorem}

The proof of the last result follows exactly the same arguments as
those of Mikami and Thieullen \cite{Mikami2006} in the proof of their
Theorem 2.1. We report it in Section~\ref{sect-theoweakdual} for
completeness.

\subsection{An enlarged space}
\label{subsecenlargedspace}

In preparation of the proof of Theorem~\ref{theodual}, we introduce
the enlarged canonical space
%
\begin{equation}
\label{eqenlargedspace} \Ob:= C \bigl([0,1], \R^d\times
\R^{d^2}\times\R^{d} \bigr),
\end{equation}
following the technique used by Haussmann \cite{Haussmann1985} in the
proof of his Proposition~3.1.

On $\Ob$, we denote the canonical filtration by $ \Fb= ( \Fcb_t)_{0\le
t \le1}$ and the canonical process by $(X,A,B)$, where $X$, $B$ are
$d$-dimensional processes and $A$ is a $d^2$-dimensional process.

We consider a probability measure $\Pb$ on $\Ob$ such that $X$ is an
$\Fb$-semimartin\-gale characterized by $(A,B)$ and, moreover, $(A,B)$ is
$\Pb$-a.s. absolutely continuous w.r.t. $t$ and $\nu_t \in U$, $d\Pb
\x
\,dt$-a.e., where $\nu=(\a,\beta)$ is defined by
%
\begin{equation}
\label{eqdnu} \a_t:= \limsup_{n \to\infty} n (A_t -
A_{t-{1}/{n}} ) \quad\mbox{and}\quad \beta_t:= \limsup_{n \to\infty} n
(B_t - B_{t-{1}/{n}} ).
\end{equation}
We also denote by $\Pcb$ the set of all probability measures $\Pb$ on
$(\Ob,\Fcb_1)$ satisfying the above conditions, and
\begin{eqnarray*}
\Pcb(\mu_0)&:=& \bigl\{\Pb\in\Pcb\dvtx\Pb\circ X_0^{-1}
= \mu_0 \bigr\}, \\
 \Pcb(\mu_0, \mu_1)&:=& \bigl\{
\Pb\in\Pcb(\mu_0)\dvtx\Pb\circ X_1^{-1} =
\mu_1 \bigr\}.
\end{eqnarray*}
Finally, we denote
\[
\overline{J}(\Pb):= \E^{\Pb}\int_0^1L(t,X,
\nu_t)\,dt.
\]

\begin{Lemma}\label{lemJbarlsc}
The function $\overline{J}$ is lower semicontinuous on $\Pcb$.
\end{Lemma}

\begin{pf} We follow the lines in Mikami \cite{Mikami2002}. By exactly the
same arguments for proving inequality (3.17) in \cite{Mikami2002},
under Assumptions~\ref{assumLconvex} and \ref
{assumLunifcontinuity}, we get
%
\begin{eqnarray}
\label{eqinqLlowersc} && \int_0^1 L(s, \xb,
\eta_s) \,ds
\nonumber
\\[-8pt]
\\[-8pt]
\nonumber
&&\qquad \ge \frac{1}{1+ \Delta_t L(\eps)} \int_0^{1-\eps} L \biggl(s,
\xb,\frac{1}{\eps}\int_s^{s+\eps}
\eta_t \,dt \biggr)\,ds - \Delta_t L (\eps)
\end{eqnarray}
for every $\eps< 1$, $\xb\in\O$ and $\R^{d^2+d}$-valued process
$\eta$.

Suppose now $(\Pb{}^n)_{n \ge1}$ is a sequence of probability measures
in $\Pcb$ which converges weakly to some $\Pb{}^0\in\Pcb$. Replacing
$(\xb, \eta)$ in \eqref{eqinqLlowersc} by $(X, \nu)$, taking expectation
under $\Pb{}^n$, by the definition of $\nu_t$ in \eqref{eqdnu} as well
as the absolute continuity of $(A,B)$ in $t$, it follows that
\begin{eqnarray*}
 \overline{J} \bigl(\Pb{}^n \bigr) &=& \E^{\Pb{}^n} \int
_0^1 L(s, X, \nu_s) \,ds
\\
&=& \frac{1}{1+\Delta_t L(\eps)} \E^{\Pb{}^n} \biggl[ \int_0^{1-\eps}
L \biggl(s,X,\frac{1}{\eps}(A_{s+\eps}-A_s),
\frac{1}{\eps}(B_{s+\eps} - B_s) \biggr) \,ds \biggr]\\
&&{} -
\Delta_t L(\eps).
\end{eqnarray*}
Next, by Fatou's lemma, we find that
\[
(X,A,B) \mapsto \int_0^{1-\eps} L \biggl(s,X,
\frac{1}{\eps}(A_{s+\eps}-A_s), \frac{1}{\eps}(B_{s+\eps}
- B_s) \biggr) \,ds
\]
is lower-semicontinuous. It follows by $\Pb{}^n \to\Pb{}^0$ that
\begin{eqnarray*}
\liminf_{n\to\infty} \overline{J} \bigl(\Pb{}^n \bigr) &\ge&
\frac{1}{1+ \Delta_t L(\eps)} \E^{\Pb{}^0} \biggl[ \int_0^{1-\eps}
L \biggl(s,X, \frac{1}{\eps}\int_s^{s+\eps}
\nu_t\,dt \biggr)\,ds \biggr] - \Delta_t L(\eps).
\end{eqnarray*}
Note that by the absolute continuity assumption of $(A,B)$ in $t$ under
$\Pb{}^0$,
\begin{eqnarray*}
\frac{1}{\eps} \int_s^{s+\eps}
\nu_t(\omega) \,dt \longrightarrow \nu_s(\omega)\qquad
\mbox{as } \eps\to0, \mbox{ for } d\Pb{}^0\x \,dt\mbox{-a.e. } (\om, s) \in\O
\x[0,1),
\end{eqnarray*}
and $\Delta_t L(\eps) \to0$ as $\eps\to0$ from Assumption \ref
{assumLunifcontinuity}; we then finish the proof by sending $\eps$
to zero and using Fatou's lemma.
\end{pf}

\begin{Remark}
In the Markovian case $L(t,\mathbf{x},u)=\ell(t, \mathbf{x}(t),u)$, for
some deterministic function $\ell$, we observe that Assumption \ref
{assumLunifcontinuity} is stronger than Assumption A2 in Mikami
\cite
{Mikami2002}. However, we can easily adapt this proof by introducing
the trajectory\vspace*{1pt} set $\{ \xb\dvtx  \sup_{0\le t,s \le1, |t-s| \le\eps}
|\xb(t) - \xb(s)| \le\delta\}$ and then letting $\eps, \delta\to
0$ as in the proof of inequality (3.17) in \cite{Mikami2002}.
\end{Remark}

Our next objective is to establish a one-to-one connection between the
cost functional $J$ defined on the set $\Pc(\mu_0, \mu_1)$ of
probability measures on $\O$ and the cost functional $\overline{J}$
defined on the corresponding set $\Pcb(\mu_0, \mu_1)$ on the enlarged
space $\Ob$.

\begin{Proposition}\label{propJJbar}
\textup{(i)} For any probability measure $\P\in\Pc(\mu_0, \mu_1)$, there
exists a probability $\Pb\in\Pcb(\mu_0, \mu_1)$ such that $J(\P)
= \Jb
(\Pb)$.\vspace*{-6pt}
\begin{longlist}[(ii)]
\item[(ii)] Conversely, let $\Pb\in\Pcb(\mu_0,\mu_1)$ be such that
$\E^{
\Pb} \int_0^1 |\beta_s| \,ds < \infty$. Then, under Assumption \ref
{assumLconvex}, there exists a probability measure $\P\in\Pc(\mu_0,\mu_1)$ such that $J(\P)\le\Jb(\Pb).$
\end{longlist}
\end{Proposition}

\begin{pf}(i) Given $\P\in\Pc(\mu_0, \mu_1)$, define the processes
$A^{\P
}$, $B^{\P}$ from decomposition \eqref{eqXdecomposition} and observe
that the mapping $\om\in\O\mapsto(X_t(\om), A_t^{\P}(\om),\break
B_t^{\P
}(\om)) \in\R^{2d+d^2}$ is measurable for every $t \in[0,1]$. Then
the mapping $\om\in\O\mapsto(X(\om),  A^{\P}(\om), B^{\P}(\om))
\in
\Ob$
is also measurable; see, for example, discussions in Chapter~2 of
Billingsley \cite{Billinsley1968} at page 57.

Let $\Pb$ be the probability measure on $(\Ob, \Fcb_1)$ induced by
$(\P, (X, A^{\P}(X),\break  B^{\P}(X)))$. In the enlarged space $(\Ob, \Fcb_1,
\Pb
)$, the canonical process $X$ is clearly a continuous semimartingale
characterized by $(A^{\P}(X), B^{\P}(X))$. Moreover, $(A^{\P}(X),
B^{\P
}(X) ) = (A,B)$, $\Pb$-a.s., where $(X,A,B)$ are canonical processes in~$\Ob$.
It follows that, on the enlarged space $(\Ob, \Fb, \Pb)$, $X$ is
a continuous semimartingale characterized by\vadjust{\goodbreak} $(A,B)$. Also, $(A,B)$ is
clearly $\Pb$-a.s. absolutely continuous in $t$, with $\nu^{\P}(X)_t =
\nu_t$, $d \Pb\x \,dt$-a.e., where $\nu$ is defined in \eqref{eqdnu}.
Then $\Pb$ is the required probability in $\Pcb(\mu_0, \mu_1)$ and
satisfies $\Jb(\Pb) = J(\P)$.

(ii) Let us first consider the enlarged space $\Ob$, and denote by
$\Fb{}^X = (\Fcb{}^X_t)_{0 \le t \le1}$ the filtration generated by process
$X$. Then for every $\Pb\in\Pcb(\mu_0, \mu_1)$, $(\Ob, \Fb{}^X,
\Pb,
X)$ is still a continuous semimartingale, by the stability property of
semimartingales. It follows from Theorem~\ref{theoIto2Diffusion} in
the \hyperref[app]{Appendix} that the decomposition of $X$ under filtration $\Fb{}^X =
(\Fcb{}^X_t)_{0 \le t \le1}$ can be written as
\[
X_t = X_0 + \bar{B}(X)_t +
\bar{M}(X)_t = X_0 + \int_0^t
\bar{\beta}_s \,ds + \bar{M}(X)_t,
\]
with $\Ab(X)_t:= \langle\bar{M}(X) \rangle_t = \int_0^t \bar{\a
}_s \,ds$,
$\bar{\beta}_s = \E^{\Pb}  [ \beta_s  | \Fcb{}^X_s
]$ and
$ \bar{\a}_s = \a_s$, $d \Pb\x \,dt$-a.e. Moreover, by the convexity property
\eqref{eqconvexitysetSigma} of the set $U$, it follows that $(\bar
\alpha,\bar{\beta}) \in U$, $d \Pb\x \,dt$-a.e. Finally, since $\Fcb{}^{X}_t = \Fc_t \otimes \{ \varnothing, C([0,1], \R^{d^2} \x\R^{d})
\}$, $\Pb$ then induces a probability measure $\P$ on $(\O, \Fc_1)$ by
\[
\P[E]:= \Pb \bigl[E \x C\bigl([0,1], \R^{d^2} \x\R^d\bigr)
\bigr] \qquad \forall E \in\Fc_1.
\]
Clearly, $\P\in\Pc(\mu_0, \mu_1)$ and $J(\P) \le\Jb(\Pb)$ by the
convexity of $L$ in $b$ of Assumption~\ref{assumLconvex} and Jensen's
inequality.
\end{pf}

\begin{Remark} \label{remEbetafini}
Let $\Pb\in\Pcb$ be such that $\Jb(\Pb) < \infty$, then from the
coercivity property of $L$ in $u$ in Assumption \ref
{assumLcroissantp}, it follows immediately that
$\E^{ \Pb} \int_0^1 |\beta_s| \,ds < \infty.$
\end{Remark}

\subsection{Lower semicontinuity and existence}

By the correspondence between $J$ and $\overline{J}$ (Proposition \ref
{propJJbar}) and the lower semicontinuity of $\overline{J}$ (Lemma
\ref{lemJbarlsc}), we now obtain the corresponding property for $V$ under
the crucial Assumption~\ref{assumLcroissantp}, which guarantees the
tightness of any minimizing sequence of our problem $V(\mu_0, \mu_1)$.

\begin{Lemma}\label{lemVlsc}
Under Assumptions~\ref{assumLconvex},~\ref{assumLunifcontinuity}
and~\ref{assumLcroissantp}, the map
\[
(\mu_0, \mu_1) \in\mathbf{ M}\bigl(\R^d\bigr)
\x\mathbf{ M}\bigl(\R^d\bigr) \longmapsto V(\mu_0,
\mu_1) \in\overline{\R}:= \R\cup\{\infty\}
\]
is lower semicontinuous.
\end{Lemma}

\begin{pf} \hspace*{-0.5pt}We follow the arguments in Lemma 3.1 of Mikami and Thieul\-len~\cite{Mikami2006}.
Let $(\mu_0^n)$ and $(\mu_1^n)$ be two sequences in
$\mathbf{ M}(\R^d)$ converging weakly to $\mu_0,\mu_1\in\mathbf{ M}(\R^d)$,
respectively, and let us prove that
\[
\liminf_{n \to\infty} V\bigl(\mu_0^n,
\mu_1^n\bigr) \ge V(\mu_0,
\mu_1).
\]
We focus on the case $\liminf_{n \to\infty} V(\mu_0^n, \mu_1^n)<\infty
$, as the result is trivial in the alternative case. Then, after
possibly extracting a subsequence, we can assume that $(V(\mu^n_0, \mu^n_1))_{n \ge1}$
is bounded, and there is a sequence\vadjust{\goodbreak} $(\P_n )_{n \ge
1}$ such that $\P_n\in\Pc(\mu_0^n, \mu_1^n)$ for all $n\ge1$ and
%
\begin{equation}
\label{VJPn} \qquad\sup_{n \ge1} J(\P_n) < \infty,\qquad  0 \le J(
\P_n) - V\bigl(\mu^n_0, \mu^n_1
\bigr) \longrightarrow0 \qquad \mbox{as } n \to\infty.
\end{equation}
By Assumption~\ref{assumLcroissantp} it follows that $\sup_{n \ge1}
\E^{\P_n} \int_0^1 |\nu_s^{\P_n}|^p \,ds < \infty$. Then, it
follows from
Theorem 3 of Zheng \cite{Zheng1985} that the sequence $(\Pb_n)_{n
\ge
1}$, of probability measures induced by $(\P_n, X, A^{\P_n}, B^{\P_n})$
on $(\Ob, \Fcb_1)$, is tight. Moreover, under any one of their limit
laws $\Pb$, the canonical process $X$ is a semimartingale characterized
by $(A,B)$ such that $(A,B)$ are still absolutely continuous in~$t$.
Moreover, $\nu\in U, d\Pb\x \,dt$-a.e. since $\frac{1}{t-s} (A_t -
A_s, B_t - B_s) \in U, d\Pb$-a.s. for every $t,s \in[0,1]$, hence,
$\Pb\in\Pcb(\mu_0, \mu_1)$. We then deduce from (\ref{VJPn}),
Proposition~\ref{propJJbar} and Lemma~\ref{lemJbarlsc} that
\[
\liminf_{n \to\infty} V\bigl(\mu_0^n,
\mu_1^n\bigr) = \liminf_{n \to\infty} J(
\P_n) = \liminf_{n \to\infty} \overline{J}(\Pb_n) \ge
\overline{J}(\Pb).
\]
By Remark~\ref{remEbetafini} and Proposition~\ref{propJJbar}, we may
find $\P\in\Pc(\mu_0, \mu_1)$ such that $\Jb(\Pb)\ge J(\P)$. Hence,
$\liminf_{n \to\infty} V(\mu_0^n, \mu_1^n)\ge J(\P)\ge V(\mu_0,\mu_1)$, completing the proof.
\end{pf}

\begin{Proposition}\label{propexisminimizer}
Let Assumptions~\ref{assumLconvex},~\ref{assumLunifcontinuity} and
\ref{assumLcroissantp} hold true. Then for every $\mu_0$, $\mu_1
\in
\mathbf{M}(\R^d)$ such that $V(\mu_0, \mu_1) < \infty$, existence holds
for the minimization problem $V(\mu_0, \mu_1)$. Moreover, the set of
minimizers $ \{ \P\in\Pc(\mu_0, \mu_1)\dvtx J(\P) = V(\mu_0,
\mu_1)
\}$ is a compact set of probability measures on $\O$.
\end{Proposition}

\begin{pf} We just let $(\mu_0^n, \mu_1^n) = (\mu_0, \mu_1)$ in the proof
of Lemma~\ref{lemVlsc}, then the required existence result is proved
by following the same arguments.
\end{pf}

\subsection{Convexity}

\begin{Lemma}\label{lemmconvex}
Let Assumptions~\ref{assumLconvex} and~\ref{assumLcroissantp}
hold, then the map $(\mu_0, \mu_1) \mapsto V(\mu_0, \mu_1)$ is convex.
\end{Lemma}

\begin{pf} Given $\mu_0^1$, $\mu_0^2$, $\mu_1^1$, $\mu_1^2 \in\mathbf{
M}(\R^d)$ and $\mu_0 = \theta\mu_0^1 + (1 - \theta) \mu_0^2$, $\mu_1 =
\theta\mu_1^1 + (1 - \theta) \mu_1^2$ with $\theta\in(0,1)$, we
shall prove that
\[
V(\mu_0, \mu_1) \le \theta V\bigl(
\mu_0^1, \mu_1^1\bigr) + (1 -
\theta) V\bigl(\mu_0^2, \mu_1^2
\bigr).
\]
It is enough to show that for both $\P_i \in\Pc(\mu_0^i, \mu_1^i)$
such that $ J(\P_i) < \infty$, $i = 1,2$, we can find $\P\in\Pc
(\mu_0, \mu_1)$ satisfying
%
\begin{equation}
\label{eqconvexityJ} J(\P)\le\theta J(\P_1) + (1-\theta) J(
\P_2).
\end{equation}
As in Lemma~\ref{lemVlsc}, let us consider the enlarged space $\Ob$,
on which the probability measures $\Pb_i$ are induced by $(\P_i, X,
A^{\P_i}, B^{\P_i})$, $i=1,2$. By Proposition~\ref{propJJbar}, $(\Pb_i)_{i=1,2}$
are probability measures under which $X$ is a $\Fb$-semimar\-tingale
characterized by the same process $(A,B)$, which is
absolutely continuous in $t$, such that $J(\P_i) = \Jb(\Pb_i), i=1,2.$

By Corollary III.2.8 of Jacod and Shiryaev \cite{Jacod1987}, $\Pb:=
\theta\Pb_1 + (1 - \theta) \Pb_2$ is also a probability measure under
which $X$ is an $\Fb$-semimartingale characterized by $(A,B)$. Clearly,
$\nu\in U, d\Pb\x \,dt$-a.e. since it is true $d \Pb_i \x \,dt$-a.e. for
$i = 1,2$. Thus, $\Pb\in\Pc(\mu_0, \mu_1)$ and it satisfies that
\begin{eqnarray*}
\Jb(\Pb) = \theta\Jb(\Pb_1) + (1-\theta) \Jb(\Pb_2) =
\theta J(\P_1) + (1 - \theta) J(\P_2) < \infty.
\end{eqnarray*}

Finally, by Remark~\ref{remEbetafini} and Proposition \ref
{propJJbar}, we can construct $\P\in\Pc(\mu_0, \mu_1)$ such that
$J(\P) \le\Jb(\Pb)$, and it follows that inequality \eqref
{eqconvexityJ} holds true.
\end{pf}

\subsection{Proof of the duality result}

We follow the first part of the proof of Theorem 2.1 in Mikami and
Thieullen \cite{Mikami2006}. If $V(\mu_0, \mu_1)$ is infinite for
every $\mu_1 \in\mathbf{M}(\R^d)$, then $J(\P) = \infty$ for all
$\P
\in\Pc(\mu_0)$. It follows from \eqref{eqVdual} and Lemma~\ref
{lemmdynamprogmu0} that
\[
V(\mu_0, \mu_1) = \Vc(\mu_0,
\mu_1) = \infty.
\]

Now, suppose that $V(\mu_0, \cdot)$ is not always infinite. Let $\overline{\mathbf{M}}(\R^d)$ be the space of all finite signed measures on $(\R^d,
\Bc(\R^d))$, equipped with weak topology, that is, the coarsest topology
making $\mu\mapsto\mu(\phi)$ continuous for every $\phi\in C_b(\R^d)$.
As indicated in Section 3.2 of \cite{Deuschel1989}, the topology
inherited by $\mathbf{ M}(\R^d)$ as a subset of $\overline{\mathbf{M}}(\R^d)$ is its
weak topology. We then extend $V(\mu_0, \cdot)$ to $\overline{\mathbf{M}}(\R^d)
\supset\mathbf{ M}(\R^d)$ by setting $V(\mu_0, \mu_1) = \infty$ when
$\mu_1 \in\overline{\mathbf{M}}(\R^d) \setminus\mathbf{M}(\R^d)$, thus, $\mu_1
\mapsto V(\mu_0, \mu_1)$ is a convex and lower semicontinuous function
defined on $\overline{\mathbf{M}}(\R^d)$. Then, the duality result $V= \Vc$ follows
from Theorem 2.2.15 and Lemma 3.2.3 in \cite{Deuschel1989}, together
with the fact that for $\lambda_1 \in C_b(\R^d)$,
\begin{eqnarray*}
&& \sup_{\mu_1 \in\mathbf{M}(\R^d)} \bigl\{ \mu_1(-\lambda_1) - V(
\mu_0, \mu_1) \bigr\}
\\
&& \qquad=- \mathop{\inf_{\mu_1 \in\mathbf{M}(\R^d)}}_{
\P\in\Pc(\mu_0, \mu_1)} \E^{\P} \biggl[ \int_0^1
L\bigl(s,X, \nu^{\P}_s\bigr)\,ds + \lambda_1(X_1)
\biggr]
\\
&&\qquad= - \inf_{\P\in\Pc(\mu_0)} \E^{\P} \biggl[ \int_0^1
L\bigl(s,X, \nu^{\P}_s\bigr)\,ds + \lambda_1(X_1)
\biggr]
\\
&&\qquad= - \mu_0 (\lambda_0),
\end{eqnarray*}
where the last equality follows by Lemma~\ref{lemmdynamprogmu0}.

\subsection{\texorpdfstring{Proof of Theorem \protect\ref{theoweakdual}}{Proof of Theorem 3.8}}
\label{sect-theoweakdual}

The proof is almost the same as that of Theorem~2.1 of Mikami and
Thieullen \cite{Mikami2006}; we report it here for completeness. Let
$\psi\in C_c^{\infty}([-1,1]^d, \R^+)$ be such that $\int_{\R^d}
\psi
(x) \,dx = 1$, and define $ \psi_{\eps}(x):= \eps^{-d}\psi(x/\eps
)$. We
claim that
%
\begin{equation}
\label{eqdualapprox} \overline{\Vc} (\mu_0, \mu_1) \ge
\frac{ \Vc(\psi_{\eps} \ast\mu_0, \psi_{\eps} \ast\mu_1)} {
1+\Delta_x L(\eps)} - \Delta_x L(\eps).
\end{equation}
Since the inequality $\Vc\ge\overline{\Vc}$ is obvious, the required
result is then obtained by sending $\eps\to0$ and using Assumption
\ref{assumLunifcontinuityinX} together with Lemma~\ref{lemVlsc}.

Hence, we only need to prove the claim \eqref{eqdualapprox}. Let us
denote $\delta:=\Delta_xL(\eps)$ in the rest of this proof. We first
observe from Assumption~\ref{assumLunifcontinuityinX} that
\begin{eqnarray*}
L(s,\xb,u) &\ge& \frac{L(s,\xb+z,u)}{1 + \delta} - \delta\qquad \mbox{for all } z \in\R
\mbox{ satisfying } |z| \le\eps,
\end{eqnarray*}
where $ \xb+z:= (\xb(t) + z)_{0 \le t \le1} \in\O$.
For an arbitrary $\lambda_1 \in C_b(\R^d)$, we denote $\lambda_1^{\eps}:= (1+\delta)^{-1}\lambda_1 \ast\psi_{\eps} \in C_b^{\infty}$.
Let $\P\in\Pc(\mu_0)$ and $Z$ be a r.v. independent of $X$ with
distribution defined by the density function $\psi_{\eps}$ under $\P$.
Then the probability $\Pb_\eps$ on $\Ob$ induced by $( \P, X+Z:= (X_t
+ Z)_{0 \le t \le1}, A^{\P}, B^{\P} )$ is in $\Pcb(\psi_{\eps}
\ast\mu_0)$, and
\begin{eqnarray*}
&& \E^{\P} \biggl[ \int_0^1 L
\bigl(s,X,\nu^{\P}_s\bigr)\,ds + \lambda_1^{\eps}(X_1)
\biggr]
\\
&&\qquad\ge - \delta + \frac{1}{1+\delta} \E^{\P} \biggl[ \int
_0^1 L\bigl(s,X+Z, \nu^{\P}_s
\bigr)\,ds + \lambda_1(X_1+Z) \biggr]
\\
&&\qquad= - \delta + \frac{1}{1+\delta} \E^{\Pb_\eps} \biggl[ \int
_0^1 L(s,X, \nu_s)\,ds +
\lambda_1(X_1) \biggr]
\\
&&\qquad \ge - \delta + \frac{1}{1+\delta} \inf_{\tilde{\P} \in\Pc(\psi_{\eps} \ast\mu_0)} \E^{\tilde{\P}}
\biggl[ \int_0^1 L\bigl(s,X,
\nu^{\tilde{\P}}_s\bigr)\,ds + \lambda_1(X_1)
\biggr],
\end{eqnarray*}
where the last inequality follows from Proposition~\ref{propJJbar}.

Notice that $\mu_1(\lambda^{\eps}_1) = (1+\delta)^{-1}(\psi_\eps
\ast
\mu_1)(\lambda_1)$ by Fubini's theorem. Then, by the arbitrariness of
$\lambda_1 \in C_b(\R^d)$ and $\P\in\Pc(\mu_0)$, the last inequality
implies~\eqref{eqdualapprox}.

\section{Characterization of the dual formulation}
\label{secPDEdynprog}

In the rest of the paper we assume that
\[
L(t, \xb, u) = \ell\bigl(t, \xb(t), u\bigr),
\]
where the deterministic function $\ell\dvtx (t,x,u) \in[0,1] \x\R^d \x U
\mapsto\ell(t,x,u) \in\R^+$ is nonnegative and convex in $u$. Then,
the function $\lambda_0$ in \eqref{eqlambda0L} is reduced to the
value function of a standard Markovian stochastic control problem:
%
\begin{equation}
\label{eqlambda0l} \lambda_0(x) = \inf_{\P\in\Pc(\delta_x)}
\E^{\P} \biggl[\int_0^1 \ell\bigl(s,
X_s, \nu^{\P}_s\bigr) \,ds +
\lambda_1(X_1) \biggr].
\end{equation}
Our main objective is to characterize $\lambda_0$ by means of the
corresponding dynamic programming equations. Then in the case of
bounded characteristics, we show more regularity as well as
approximation properties of $\lambda_0$, which serves as a preparation
for the numerical approximation in Section~\ref{secnumericalapprox}.

\subsection{PDE characterization of the dynamic value function}

Let us consider the probability measures $\P$ on the canonical space
$(\O, \Fc_1)$, under which the canonical process $X$ is a
semimartingale on $[t,1]$, characterized by $\int_t^{\cdot} \nu^{\P}_s
\,ds$ for some progressively measurable process $\nu^{\P}$. As discussed
in Remark~\ref{remuniquenu}, $\nu^{\P}$ is unique on $\O\x[t,1]$ in
the sense of $d \P\x \,dt$-a.e.
Following the definition of $\Pc$ just below~\eqref
{eqconvexitysetSigma}, we denote by $\Pc_t$ the collection of all
such probability measures $\P$ such that $\nu^{\P}_s \in U$, $d \P
\x
\,dt$-a.e. on $\O\x[t,1]$. Let
%
\begin{equation}
\Pc_{t,x}:= \bigl\{ \P\in\Pc_t \dvtx
\P[X_s = x, 0 \le s \le t] = 1 \bigr\}.
\end{equation}
We notice that under probability $\P\in\Pc_{t,x}$, $X$ is a
semimartingale with $\nu^{\P}_s = 0$, $d\P\x \,dt$-a.e. on $\O\x
[0,t]$. The dynamic value function is defined for any $\lambda_1\in
C_b(\R^d)$ by
%
\begin{equation}
\label{eqlambdal} \lambda(t,x):= \inf_{\P\in\Pc_{t,x}} \E^{\P}
\biggl[ \int_t^1 \ell\bigl(s,X_s,
\nu^{\P}_s\bigr) \,ds + \lambda_1(X_1)
\biggr].
\end{equation}
As in the previous sections, we also introduce the corresponding
probability measures on the enlarged space $(\Ob, \Fcb_1)$. For all $t
\in[0,1]$, we denote by $\Pcb_t$ the collection of all probability
measures $\Pb$ on $(\Ob, \Fcb_1)$ under which $X$ is a semimartingale
characterized by $(A,B)$ in $\Ob$ and $\nu\in U$, $d\Pb\x \,dt$-a.e. on
$\O\x[t,1]$, where $\nu$ is defined above \eqref{eqdnu}. For every
$(t,x,a,b) \in[0,1] \x\R^d \x\R^{d^2} \x\R^d$, let
%
\begin{eqnarray}
\Pcb_{t,x,a,b}:= \bigl\{ \Pb\in\Pcb\dvtx \Pb \bigl[(X_s,
A_s, B_s) = (x,a,b), 0 \le s \le t \bigr]=1 \bigr\}.
\end{eqnarray}
By similar arguments as in Proposition~\ref{propJJbar}, we have under
Assumption~\ref{assumLconvex} that
%
\begin{equation}
\label{eqlambdalrelax} \lambda(t,x) = \inf_{\Pb\in\Pcb_{t,x,a,b}} \E^{\Pb}
\biggl[ \int_t^1 \ell(s,X_s,
\nu_s) \,ds + \lambda_1(X_1) \biggr]
\end{equation}
for all $(a,b) \in\R^{d^2} \x\R^d$.

We would like to characterize the dynamic value function $\lambda$ as
the viscosity solution of a dynamic programming equation. The first
step is as usual to establish the dynamic programming principle (DPP).
We observe that a weak dynamic programming principle as introduced in
Bouchard and Touzi~\cite{Bouchard2011} suffices to prove that
$\lambda
$ is a viscosity solution of the corresponding dynamic programming
equation. The main argument in \cite{Bouchard2011} to establish the
weak DPP is the conditioning and pasting techniques of the control
process, which is convenient to use for control problems in a strong
formulation, that is, when the measure space $(\O, \Fc)$ as well as the
probability measure $\P$ are fixed a priori. However, we cannot use
their techniques since our problem is in weak formulation, where the
controlled process is fixed as a canonical process and the controls are
given as probability measures on the canonical space.\looseness=-1

We will prove the standard dynamic programming principle. For a simpler
problem (bounded convex controls\vadjust{\goodbreak} set $U$ and bounded cost functions,
etc.), a DPP is shown (implicitly) in Haussmann \cite{Haussmann1985}.
El Karoui, Nguyen and JeanBlanc \cite{ElKaroui1987} considered a
relaxed optimal control problem and provided a scheme of proof without
all details. Our approach is to adapt the idea of~\cite{ElKaroui1987}
in our context and to provide all details for their scheme of proof.

\begin{Proposition}\label{propdpp}
Let Assumptions~\ref{assumLconvex},~\ref{assumLunifcontinuity},
\ref{assumLcroissantp} hold true. Then, for all $\Fb$-stopping time
$\tau$ with values in $[t,1]$, and all $(a,b) \in\R^{d^2+d}$,
\[
\lambda(t,x) = \inf_{\Pb\in\Pcb_{t,x,a,b}} \E^{\P} \biggl[ \int
_t^{\tau} \ell(s,X_s,
\nu_s) \,ds + \lambda(\tau, X_{\tau}) \biggr].
\]
\end{Proposition}

The proof is reported in Section~\ref{sect-proofDPP}.
The dynamic programming equation is the infinitesimal version of the
above dynamic programming principle.
Let
%
\begin{equation}
\label{eqHamiltonian} H(t,x,p,\Gamma):= \inf_{ (a,b) \in U} \biggl[ b \cdot p
+ \frac{1}{2} a \cdot\Gamma+ \ell(t,x,a,b) \biggr]
\end{equation}
for all $(p,\Gamma) \in\R^d \x S_d.$

\begin{Theorem}\label{theouviscositysol}
Let Assumptions~\ref{assumLconvex},~\ref{assumLunifcontinuity},
\ref{assumLcroissantp} hold true, and assume further that $\lambda$
is locally bounded and $H$ is continuous. Then, $\lambda$ is a
viscosity solution of the dynamic programming equation
%
\begin{equation}
\label{eqHJB} - \partial_t \lambda(t,x) - H\bigl(t,x, D \lambda,
D^2 \lambda\bigr) =0,
\end{equation}
with terminal condition $\lambda(1,x) = \lambda_1(x)$.
\end{Theorem}

The proof is very similar to that of Corollary 5.1 in \cite
{Bouchard2011}; we report it in the \hyperref[app]{Appendix} for completeness.

\begin{Remark}
We first observe that $H$ is concave in $(p, \Gamma)$ as infimum of a
family of affine functions. Moreover, under Assumption \ref
{assumLcroissantp}, $\ell$ is positive and $u \mapsto\ell(t,x,u)$
has growth larger than $|u|^p$ for $p > 1$; it follows that $H$ is
finite valued and hence continuous in $(p,\Gamma)$ for every fixed
$(t,x) \in[0,1] \x\R^d$. If we assume further that $(t,x) \mapsto
\ell
(t,x,u)$ is uniformly continuous uniformly in $u$, then clearly $H$ is
continuous in $(t,x,p,\Gamma)$.
\end{Remark}

\begin{Remark}
The following are two sets of sufficient conditions to ensure the local
boundedness of $\lambda$ in \eqref{eqlambdal}.
\begin{longlist}[(ii)]
\item[(i)] Suppose $0 \in U$, and let $\P\in\Pc_t$ be such that $\nu^{\P}_s
= 0$, $d\P\x \,dt$-a.e. Then, $\lambda(t,x) \le|\lambda_1|_{\infty} +
\int_t^1 \ell(s,x,0) \,ds$ and, hence, $\lambda$ is locally bounded.

\item[(ii)] Suppose that there are constants $C>0$ and $(a_0, b_0) \in U$ such
that $\ell(t,x,a_0,b_0) \le C e^{C |x|}$, for all $(t,x) \in[0,1] \x
\R^d$. By considering $\P\in\Pc_t$ induced by the process $Y = (Y_s)_{t
\le s \le1}$ with $Y_s:= x + b_0(s-t) + a_0^{1/2} (W_s - W_t)$, it
follows that $\lambda(t,x) \le|\lambda_1|_{\infty} + \E[C e^{C
\max
_{t \le s \le1} |Y_s| }] < \infty$.
\end{longlist}
\end{Remark}

\subsubsection{Proof of the dynamic programming principle}
\label{sect-proofDPP}

We first prove that the dynamic value function $\lambda$ is measurable
and we can choose ``in a measurable way'' a family of probabilities
$(\Q_{t,x,a,b})_{(t,x,a,b) \in[0,1] \x\R^{2d+d^2}}$ which achieves
(or achieves with $\eps$ error) the infimum in \eqref
{eqlambdalrelax}. The main argument is Theorem \ref
{TheoMeasurabilitySup} cited in the \hyperref[app]{Append}ix which follows directly
from the measurable selection theorem.

Let $\lambda^*$ be the upper semicontinuous envelope of the function
$\lambda$, and
\begin{eqnarray*}
\tilde{\Pc}_{t,x,a,b} &:=& \biggl\{\Pb\in\Pcb_{t,x,a,b} \dvtx
\E^{\Pb} \biggl[ \int_t^1 \ell
(s,X_s,\nu_s) \,ds + \lambda_1(X_1)
\biggr] \le\lambda^*(t,x) \biggr\},
\\
 \tilde{\Pc} &:=& \bigl\{ (t,x,a,b, \Pb) \dvtx\Pb\in\tilde{\Pc}_{t,x,a,b}
\bigr\}.
\end{eqnarray*}
In the following statement, for the Borel $\sigma$-field $ \Bc
([0,1]\x
\R^{2d+d^2})$ of $[0,1] \x\R^{2d+d^2}$ with an arbitrary probability
measure $\mu$ on it, we denote by $\Bc^{\mu}([0,1]\x\R^{2d+d^2})$ its
$\sigma$-field completed by $\mu$.

\begin{Lemma}\label{lemmmeasurableselection}
Let Assumptions~\ref{assumLconvex},~\ref{assumLunifcontinuity},
\ref{assumLcroissantp} hold true, and assume that $\lambda$ is
locally bounded. Then, for any probability measure $\mu$ on $ (
[0,1] \x\R^{2d+d^2},\break \Bc([0,1]\x\R^{2d+d^2}) )$,
\begin{longlist}[(ii)]
\item[(i)] the function $(t,x,a,b) \mapsto\lambda(t,x)$ is $\Bc^{\mu
}([0,1]\x\R^{2d+d^2})$-measurable,

\item[(ii)] for any $\eps>0$, there is a family of probability $(\Qb{}^{\eps}_{t,x,a,b})_{(t,x,a,b) \in[0,1] \x\R^{2d + d^2}}$ in $\tilde
{\Pc}$ such that $(t,x,a,b) \mapsto\Qb{}^{\eps}_{t,x,a,b}$ is a
measurable map from $[0,1] \times\R^{2d+d^2}$ to $\mathbf{ M}(\Ob)$ and
\begin{eqnarray*}
\E^{\Qb{}^{\eps}_{t,x,a,b}} \biggl[ \int_t^1
\ell(s,X_s,\nu_s) \,ds + \lambda_1(X_1)
\biggr] &\le& \lambda(t,x) + \eps, \qquad\mu\mbox{-a.s.}
\end{eqnarray*}
\end{longlist}
\end{Lemma}

\begin{pf} By Lemma~\ref{lemJbarlsc}, the map $\Pb\mapsto\E^{\Pb}
[\int_t^1 \ell(s,X_s,\nu_s) \,ds + \lambda_1(X_1)  ]$ is lower
semicontinuous, and therefore measurable. Moreover, $ \tilde{\Pc
}_{t,x,a,b}$ is non\-empty for every $(t,x,a,b) \in[0,1] \x\R^{2d+d^2}$. Finally, by using the same arguments as in the proof of
Lemma~\ref{lemVlsc}, we see that $\tilde{\Pc}$ is a closed subset of
$[0,1] \x\R^{2d + d^2} \x\mathbf{ M}(\Ob)$. Then, both items of the lemma
follow from Theorem~\ref{TheoMeasurabilitySup}.
\end{pf}

We next prove the stability properties of probability measures under
conditioning and concatenations at stopping times, which will be the
key-ingredients for the proof of the dynamic programming principle.

We first recall some results from Stroock and Varadhan \cite
{Stroock1979} and define some notation:
\begin{itemize}
\item For $ 0 \le t \le1$, let $\Fcb_{t,1}:= \sigma( (X_s,A_s,B_s)
\dvtx t\le s\le1)$, and let $\Pb$ be a probability measure on $(\Ob,
\Fcb_{t,1})$ with $\Pb[(X_t,A_t,B_t) = \eta_t] = 1$ for some $\eta\in
C([0,t], \R^{2d + d^2})$. Then, there is a unique probability measure
$\delta_{\eta} \otimes_t \Pb$ on $(\Ob, \Fcb_1)$ such that
$\delta_{\eta
} \otimes_t \Pb[(X_s,A_s,B_s) = \eta_s,0 \le s \le t] = 1$ and
$\delta_{\eta} \otimes_t \Pb[A] = \Pb[A]$ for all $A \in\Fcb_{t,1} $. In
addition, if $\Pb$ is also a probability measure on $(\Ob, \Fcb_1)$,
under which a process $M$ defined on $\Ob$ is a $\Fb$-martingale after
time $t$, then $M$ is still a $\Fb$-martingale after time $t$ in
probability space $(\Ob, \Fcb_1, \eta\otimes_t \Pb)$. In particular,
for $t\in[0,1]$, a constant $c_0\in\R^{2d + d^2}$ and $\Pb$ satisfying
$\Pb[(X_t,A_t,B_t) = c_0] = 1$, we denote $\delta_{c_0} \otimes_t
\Pb:=\delta_{\eta^{c_0}} \otimes_t \Pb$, where $\eta^{c_0}_s=c_0,
s\in[0,t]$.
\item Let $\Qb$ be a probability measure on $(\Ob, \Fcb_1)$ and
$\tau$
a $\Fb$-stopping time. Then, there is a family of probability measures
$(\Qb_{\om})_{\om\in\Ob}$ such that $\om\mapsto\Qb_{\om}$ is
$\Fcb_{\tau}$-measurable, for every $E \in\Fcb_1$, $\Qb[E | \Fcb_{\tau
}](\om
) = \Qb_{\om}[E]$ for $\Qb$-almost every $\om\in\Ob$ and, finally,
$\Qb_{\om}[(X_t,A_t,B_t) = \om_t \dvtx t\le\tau(\om)] = 1$, for all
$\om
\in\Ob$. This is Theorem 1.3.4 of \cite{Stroock1979}, and $(\Qb_{\om
})_{\om\in\Ob}$ is called the regular conditional probability
distribution (r.c.p.d.)
\end{itemize}

\begin{Lemma}\label{lemmrcpdsemimartingale}
Let $\Pb\in\Pcb_{t,x,a,b}$, $\tau$ be an $\Fb$-stopping time taking
value in $[t,1]$, and $(\Qb_{\om})_{\om\in\Ob}$ be a r.c.p.d. of
$\Pb
| \Fcb_{\tau}$. Then there is a $\Pb$-null set $N \in\Fcb_{\tau}$ such
that $\delta_{\om_{\tau(\om)}} \otimes_{\tau(\om)} \Qb_{\om}
\in\Pcb_{\tau(\om), \om_{\tau(\om)}}$ for all $\om\notin N$.
\end{Lemma}

\begin{pf} Since $\Pb\in\Pcb_{t,x,a,b}$, it follows from Theorem II.2.21
of Jacod and Shiryaev \cite{Jacod1987} that
\[
(X_s - B_s)_{t \le s \le1}, \qquad\bigl((X_s
- B_s)^2 - A_s \bigr)_{t \le
s \le1}
\]
are all local martingales after time $t$. Then it follows from Theorem
1.2.10 of Stroock and Varadhan \cite{Stroock1979} together with a
localization technique that there is a $\Pb$-null set $N_1 \in\Fcb_{\tau}$ such that they are still local martingales after time $\tau
(\om
)$ both under $\Qb_{\om}$ and $\delta_{\om_{\tau(\om)}} \otimes_{\tau
(\om)} \Qb_{\om}$, for all $\om\notin N_1$. It is clear, moreover,
that $\nu\in U, d\Qb_{\om} \x \,dt$-a.e. on $\Ob\x[\tau(\om), 1]$ for
$\Pb$-a.e. $\om\in\Ob$. Then there is a $\Pb$-null set $N \in\Fcb_{\tau}$ such that $\delta_{\om_{\tau(\om)}} \otimes_{\tau(\om
)} \Qb_{\om} \in\Pcb_{\tau(\om), \om_{\tau(\om)}}$ for every $\om
\notin N$.
\end{pf}

\begin{Lemma}\label{lemmconcatenation}
Let Assumptions~\ref{assumLconvex},~\ref{assumLunifcontinuity},
\ref{assumLcroissantp} hold true, and assume that $\lambda$ is
locally bounded. Let $\Pb\in\Pcb_{t,x,a,b}$, $\tau\ge t$ a $\Fb
$-stopping time, and $(\Qb_{\om})_{ \om\in\Ob}$ a family of
probability measures such that $\Qb_{\om} \in\Pcb_{\tau(\om), \om
_{\tau
(\om)}}$ and $\om\mapsto\Qb_{\om}$ is $\Fcb_{\tau}$-measurable. Then
there is a unique probability measure, denoted by $\Pb\otimes_{\tau
(\cdot)} \Qb_{\cdot}$, in $\Pcb_{t,x,a,b}$, such that $\Pb\otimes_{\tau
(\cdot)} \Qb_{\cdot} = \Pb\mbox{ on } \Fcb_{\tau}$, and
%
\begin{eqnarray}
\label{eqproprcpd} (\delta_{\om} \otimes_{\tau(\om)}
\Qb_{\om})_{\om\in\Ob}\qquad \mbox{is a r.c.p.d. of } \Pb
\otimes_{\tau(\cdot)} \Qb_{\cdot} | \Fcb_{\tau}.
\end{eqnarray}
\end{Lemma}

\begin{pf} The existence and uniqueness of the probability measure $
\Pb
\otimes_{\tau(\cdot)} \Qb_{\cdot}$ on $(\Ob, \Fcb_1)$,
satisfying \eqref
{eqproprcpd}, follows from Theorem 6.1.2 of \cite{Stroock1979}.
It remains to prove that $ \Pb\otimes_{\tau(\cdot)} \Qb_{\cdot}
\in
\Pcb_{t,x,a,b}$.

Since $\Qb_{\om} \in\Pcb_{\tau(\om),\om_{\tau(\om)}}$, $X$ is
a $\delta_{\om} \otimes_{\tau(\om)} \Qb_{\om}$-semimartingale after time
$\tau
(\om)$, characterized by $(A,B)$. Then, the processes $X-B$ and
$(X-B)^2-A$ are local martingales under $\delta_{\om} \otimes_{\tau
(\om
)} \Qb_{\om}$ after time $\tau(\om)$. By Theorem 1.2.10 of~\cite
{Stroock1979} together with a localization argument, they are still
local martingales under $\Pb\otimes_{\tau(\cdot)} \Qb_{\cdot}$. Hence,
the required result follows from Theorem II.2.21 of~\cite{Jacod1987}.
\end{pf}

We have now collected all the ingredients for the proof of the dynamic
programming principle.

\begin{pf*}{Proof of Proposition~\ref{propdpp}} Let $\tau$ be an
$\Fb
$-stopping time taking value in $[t,1]$. We proceed in two steps:
\begin{longlist}[(1)]
\item[(1)] For $\Pb\in\Pcb_{t,x,a,b}$, we denote by $(\Qb_\omega
)_{\omega
\in\Ob}$ a family of regular conditional probability distribution of
$\Pb|\Fcb_\tau$, and $\Pb_\tau^\omega:=\delta_{\om_{\tau(\om
)}} \otimes_{\tau(\om)} \Qb_{\om}$.
By the representation~\eqref{eqlambdalrelax} of $\lambda$, together
with the tower property of conditional expectations, we see that
%
\begin{eqnarray}
\label{ineq1-dpp} && \lambda(t,x)
\nonumber\\
&&\qquad= \inf_{\Pb\in\Pcb_{t,x,a,b}} \E^{\Pb} \biggl[ \int_t^{\tau}
\ell(s,X_s,\nu_s) \,ds + \int_{\tau}^1
\ell(s,X_s,\nu_s) \,ds + \lambda_1(X_1)
\biggr]
\nonumber
\\[-8pt]
\\[-8pt]
\nonumber
&&\qquad= \inf_{\Pb\in\Pcb_{t,x,a,b}} \E^{\Pb} \biggl[ \int_t^{\tau}
\ell(s,X_s,\nu_s) \,ds + \E^{\Pb_\tau^\omega} \biggl\{
\int_{\tau}^1 \ell(s,X_s,
\nu_s) \,ds + \lambda_1(X_1) \biggr\} \biggr]\hspace*{-25pt}
\\
&&\qquad\ge \inf_{\Pb\in\Pcb_{t,x,a,b}} \E^{\Pb} \biggl[ \int_t^{\tau}
\ell(s,X_s,\nu_s) \,ds + \lambda (\tau,
X_{\tau}) \biggr],
\nonumber
\end{eqnarray}
where the last inequality follows from the fact that $\Pb_\tau^\omega
\in\Pcb_{\tau(\om), \om_{\tau(\om)}}$ by Lemma~\ref
{lemmrcpdsemimartingale}.

\item[(2)] For $\eps>0$, let $(\Qb{}^{\eps}_{t,x,a,b})_{[0,1] \x\R
^{2d+d^2}}$ be the family defined in Lemma \ref
{lemmmeasurableselection}, and denote $\Qb{}^{\eps}_{\om}:= \Qb{}^{\eps
}_{\tau(\om), \om_{\tau(\om)}}$. Then $\om\mapsto\Qb{}^{\eps
}_{\om}$ is
$\Fcb_{\tau}$-measurable. Moreover, for all $\Pb\in\Pcb_{t,x,a,b}$,
we may construct by Lemmas~\ref{lemmmeasurableselection} and \ref
{lemmconcatenation} $\Pb\otimes_{\tau(\cdot)} \Qb_{\cdot} \in
\Pcb_{t,x,a,b}$ such that
\begin{eqnarray*}
&& \E^{\Pb\otimes_{\tau(\cdot)} \Qb_{\cdot}} \biggl[ \int_t^1 \ell
(s,X_s,\nu_s) \,ds + \lambda_1(
X_1) \biggr]
\\
&&\qquad \le \E^{\Pb} \biggl[ \int_t^{\tau}
\ell(s,X_s,\nu_s) \,ds + \lambda(\tau,
X_{\tau}) \biggr] + \eps.
\end{eqnarray*}
By the arbitrariness of $\Pb\in\Pcb_{t,x,a,b}$ and $\eps>0$,
together with the representation~\eqref{eqlambdalrelax} of $\lambda
$, this implies that the reverse inequality to (\ref{ineq1-dpp}) holds
true, and the proof is complete.\quad\qed
\end{longlist}
\noqed\end{pf*}

We conclude this section by the following:

\begin{pf*}{Proof of Lemma~\ref{lemmdynamprogmu0}} By the
same arguments as in Lemma~\ref{lemmmeasurableselection}, we can
easily deduce that $\lambda_0$ is $\Bc^{\mu_0}(\R^d)$-measurable, and
we just need to prove that
\[
\mu_0(\lambda_0) = \inf_{\Pb\in\Pcb(\mu_0)}
\E^{\Pb} \biggl[ \int_0^1
\ell(s,X_s, \nu_s)\,ds + \lambda_1(X_1)
\biggr].
\]
Given a probability measure $ \Pb\in\Pcb(\mu_0)$, we can get a
family of conditional probabilities $(\Qb_{\om})_{\om\in\O}$ such
that $\Qb_{\om} \in\Pcb_{0, \om_0}$, which implies that
\[
\E^{\Pb} \biggl[ \int_0^1
\ell(s,X_s, \nu_s)\,ds + \lambda_1(X_1)
\biggr] \ge \mu_0(\lambda_0) \qquad\forall\Pb\in\Pcb(
\mu_0).
\]
On the other hand, for every $\eps> 0$ and $\mu_0 \in\mathbf{ M}(\R^d)$,
we can select a measurable family of $( \Qb{}^{\eps}_x \in\Pcb_{0, x,
0,0})_{x \in\R^d}$ such that
\[
\E^{\Qb{}^{\eps}_x} \biggl[ \int_0^1
\ell(s,X_s, \nu_s)\,ds + \lambda_1(X_1)
\biggr] \le \lambda_0 (x) + \eps,\qquad \mu_0\mbox{-a.s.},
\]
and then construct a probability measure $ \mu_0 \otimes_0 \Qb{}^{\eps
}_{\cdot} \in\Pcb(\mu_0)$ by concatenation such that
\begin{eqnarray*}
\E^{ \mu_0 \otimes_0 \Qb{}^{\eps}_{\cdot}} \biggl[ \int_0^1
\ell(s,X_s, \nu_s)\,ds + \lambda_1(X_1)
\biggr] &\le& \mu_0(\lambda_0) + \eps \qquad \forall\eps> 0,
\end{eqnarray*}
which completes the proof.
\end{pf*}

\subsection{Bounded domain approximation under bounded characteristics}

The main purpose of this section is to show that when $U$ is bounded,
then $\lambda_0$ in \eqref{eqlambda0l} is Lipschitz, and we may
construct a convenient approximation of $\lambda_0$ by restricting the
space domain to bounded domains. These properties induce a first
approximation for the minimum transportation cost $V(\mu_0, \mu_1)$,
which serves as a preparation for the numerical approximation in
Section~\ref{secnumericalapprox}. Let us assume the following conditions.

\begin{Assumption}\label{assumLbounded}
The control set $U$ is compact, and $\ell$ is Lipschitz-continuous in
$x$ uniformly in $(t,u)$.
\end{Assumption}

\begin{Assumption}\label{assummu01integrability}
$ \int_{\R^d} |x| (\mu_0+\mu_1)(dx) < \infty.$
\end{Assumption}

\begin{Remark}
We suppose that $U$ is compact for two main reasons. First, the
uniqueness of viscosity solution of the HJB \eqref{eqHJB} relies on
the comparison principle, for which the boundedness of $U$ is generally
necessary. Further, to construct a convergent (monotone) numerical
scheme for a stochastic control problem, it is also generally necessary
to suppose that the diffusion functions are bounded (see also Section
\ref{subsecnumerPDE} for more discussions).
\end{Remark}

\subsubsection{The unconstrained control problem in the bounded domain}

Denote
%
\begin{eqnarray}
\label{eqdM} M &:=& \sup_{(t,x,u)\in[0,1]\x\R^d\x U} \bigl( |u| + \bigl|\ell(t,0,u)\bigr| +\bigl |
\nabla_x\ell(t,x,u)\bigr| \bigr),
\end{eqnarray}
where $\nabla_x\ell(t,x,u)$ is the gradient of $\ell$ with respect
to $x$ which exists a.e. under Assumption~\ref{assumLbounded}. Let
$O_R:= (-R, R)^d \subset\R^d$ for every $R>0$, a stopping time $\tau_R$ can be defined as the first exit time of the canonical process $X$
from~$O_R$,
\[
\tau_R:= \inf\{ t \dvtx X_t \notin O_R
\},
\]
and define for all bounded functions $\lambda_1 \in C_b(\R^d)$,
%
\begin{equation}
\label{eqlambdaR} \lambda^R(t,x):=\inf_{\P\in\Pc_{t,x}}
\E^{\P} \biggl[ \int_t^{\tau_R \wedge1} \ell
\bigl(s,X_s,\nu^{\P}_s\bigr) \,ds +
\lambda_1(X_{\tau_R\wedge1}) \biggr].
\end{equation}

\begin{Lemma}\label{lemmlambdaLip}
Suppose that $\lambda_1$ is $K$-Lipschitz satisfying $\lambda_1(0) =
0$ and Assumption~\ref{assumLbounded} holds true. Then $\lambda$ and
$\lambda^R$ are Lipschitz-continuous, and there is a constant $C$
depending on $M$ such that
\[
\bigl|\lambda(t,0)\bigr|+\bigl|\lambda^R(t,0)\bigr| + \bigl|\nabla_x
\lambda(t,x)\bigr|+\bigl|\nabla_x\lambda^R(t,x)\bigr| \le C(1 + K)
\]
for all $(t,x)\in[0,1]\times\R^d.$
\end{Lemma}

\begin{pf} We only provide the estimates for $\lambda$; those for
$\lambda_R$ follow from the same arguments. First, by Assumption \ref
{assumLbounded} together with the fact that $\lambda_1$ is
$K$-Lipschitz and $\lambda_1(0) = 0$, for every $\P\in\Pc_{t,0}$,
\begin{eqnarray*}
\E^{\P} \biggl[ \int_t^1 \ell
\bigl(s,X_s,\nu^{\P}_s\bigr) \,ds +
\lambda_1(X_1) \biggr] &\le& M + (M+K)
\sup_{t \le s \le1} \E^{\P} | X_s |.
\end{eqnarray*}
Recall that $X$ is a continuous semimartingale under $\P$ whose finite
variation part and quadratic variation of the martingale part are both
bounded by a constant $M$. Separating the two parts and using
Cauchy--Schwarz's inequality, it follows that $\E^{\P}  | X_s  |
\le M + \sqrt{M}, \forall t \le s \le1$, and then $ \llvert  \lambda
(t,0) \rrvert  \le M + (M+K)(M + \sqrt{M})$.

We next prove that $\lambda$ is Lipschitz and provide the corresponding
estimate. Observe that $\Pc_{t,y} =  \{ \P:= \tilde{\P} \circ(
X+y-x)^{-1} \dvtx\tilde{\P} \in\Pc_{t,x}  \}$. Then
\begin{eqnarray*}
&& \bigl| \lambda(t,x) - \lambda(t,y)\bigr |
\\
&&\qquad\le \sup_{\P\in\Pc_{t,x}} \E^{\P} \biggl| \int_t^1
\ell\bigl(s,X_s,\nu^{\P}_s\bigr) - \ell
\bigl(s,X_s +y-x,\nu^{\P}_s\bigr) \,ds
\\
&&\hspace*{99pt}\qquad\quad{} + \lambda_1(X_1) - \lambda_1(X_1+y-x)
\biggr|
\\
&&\qquad\le (M+K) |y-x|
\end{eqnarray*}
by the Lipschitz property of $\ell$ and $\lambda$ in $x$.\vadjust{\goodbreak}
\end{pf}

Denoting $\lambda^R_0:=\lambda^R(0,\cdot)$, in the special case where $U$
is a singleton, equation~\eqref{eqHJBR} degenerates to the heat
equation. Barles, Daher and Romano \cite{Barles1995} proved that the
error $\lambda-\lambda^R$ satisfies a large deviation estimate as $R
\to\infty$. The next result extends this estimate to our context.

\begin{Lemma} \label{lemmerrorlambdalambdaR}
Letting Assumption~\ref{assumLbounded} hold true, we denote
$|x|:=\break
\max_{i=1}^d |x_i|$ for $x \in\R^d$ and choose $R>2M$. Then, there is
a constant $C$ such that for all $K>0$, all $K$-Lipschitz function
$\lambda_1$ and $|x| \le R-M$,
\begin{eqnarray*}
\bigl|\lambda^R - \lambda\bigr|(t,x) \le C(1+K) e^{-(R-M-|x|)^2/2M}.
\end{eqnarray*}
\end{Lemma}

\begin{pf}(1) For arbitrary $(t,x)\in[0,1]\times\R^d$ and $\P
\in
\Pc_{t,x}$, we denote $Y^i:= \sup_{0 \le s \le1} \llvert  X_s^i \rrvert $,
where $X^i$ is the $i$th component of the canonical process~$X$. By the
Dubins--Schwarz time-change theorem (see, e.g., Theorem 4.6, Chapter 3
of Karatzas and Shreve \cite{Karatzas1991}), we may represent the
continuous local martingale part of $X^i$ as a time-changed Brownian
motion $W$. Since the characteristics of $X$ are bounded by $M$, we see that
%
\begin{eqnarray}
\label{eqDistributionFunctionF} S^i(R):= \P\bigl[Y^i
\ge R\bigr] &\le& \P \Bigl[ \sup_{0 \le t \le M } |W_t| \ge R -
|x_i| - M \Bigr]
\nonumber
\\
&\le& 2 \P \Bigl[\sup_{0 \le t \le M } W_t \ge R - |x_i|
- M \Bigr]
\\
&=& 4 \bigl( 1 - \mathbf{ N} \bigl(R_{|x_i|}^M \bigr) \bigr),\nonumber
\end{eqnarray}
where $R_{|x_i|}^M:=(R-M-|x_i|)/\sqrt{M}$, $\mathbf{ N}$ is the cumulative
distribution function of the standard normal distribution $N(0,1)$, and
the last equality follows from the reflection principle of the Brownian
motion. Then by integration by parts as well as \eqref
{eqDistributionFunctionF},
\begin{eqnarray*}
\label{eqCVaR} \E^{\P} \bigl[ Y^i\1_{Y^i \ge R} \bigr]
&=& R S^i(R) + \int_R^{\infty}
S^i(z) \,dz
\\
&\le& 4 \int_R^{\infty} \frac{1}{\sqrt{M}}
\frac{1}{\sqrt{2 \pi}} \exp \biggl(\frac{(z-M-|x_i|)^2}{2M} \biggr) z \,dz
\\
&=& 4\bigl(|x_i|+M\bigr) \bigl(1 - \mathbf{ N} \bigl(R_{|x_i|}^M
\bigr) \bigr) + \frac{4 \sqrt{M}}{\sqrt{2 \pi}} \exp \biggl(-\frac{(R_{|x_i|}^M)^2}{2}
\biggr).
\end{eqnarray*}
We further remark that for any $R > 0$,
\begin{eqnarray*}
\bigl(1 - \mathbf{ N}(R)\bigr) = \int_R^{\infty}
\frac{1}{\sqrt{2 \pi}} e^{
-{t^2}/{2} } \,dt \le\frac{1}{R} \int
_R^{\infty} \frac{1}{\sqrt{2 \pi
}} t e^{ -{t^2}/{2} } \,dt
= \frac{1}{\sqrt{2\pi}} \frac{1}{R} e^{{-R^2}/{2}}.
\end{eqnarray*}
(2) By definitions of $\lambda$, $\lambda^R$, it follows that for
all $(t,x)$ such that $|x| \le R-M$,
\begin{eqnarray}
\label{lambda-lambdaR} \bigl|\lambda-\lambda^R\bigr|(t,x) &\le&
\sup_{\P\in\Pc_{t,x}} \E^{\P} \biggl[ \int_{\tau_R \wedge1}^1
\bigl| \ell\bigl(s, X_s, \nu_s^{\P
}\bigr) \bigr| \,ds +
\bigl| \lambda_1(X_{\tau_R \wedge1}) - \lambda_1(X_1)
\bigr| \biggr]
\nonumber
\\
&\le& \sup_{\P\in\Pc_{t,x}} \E^{\P} \Bigl[ \Bigl( M + \sqrt{d} K R +
(M+K)\sup_{t \le s \le1}|X_s| \Bigr) \1_{\tau_R<1} \Bigr]
\nonumber
\\[-8pt]
\\[-8pt]
\nonumber
&\le& \sup_{\P\in\Pc_{t,x}} \E^{\P} \Biggl[ \sum
_{i=1}^d \bigl(M+\sqrt{d} K R+\sqrt
{d}(M+K)Y_i \bigr) \1_{Y_i \ge R} \Biggr]
\nonumber
\\
&\le& C(1+K) e^{-(R_{|x|}^M)^2/2}\nonumber
\end{eqnarray}
for some constant $C$ depending on $M$ and $d$. This completes the
proof.
\end{pf}

With the estimate in Lemma~\ref{lemmlambdaLip}, we have the
following result.

\begin{Theorem}
Suppose that Assumptions~\ref{assumLconvex}, \ref
{assumLunifcontinuity},~\ref{assumLcroissantp} hold true and $H$
given by \eqref{eqHamiltonian} is continuous. Then the function
$\lambda^R$ in \eqref{eqlambdaR} is the unique viscosity solution
of equation
%
\begin{equation}
\label{eqHJBR}\qquad - \partial_t \lambda^R(t,x) - H\bigl(t,x,
D \lambda^R, D^2 \lambda^R\bigr) = 0,\qquad
(t,x) \in[0,1) \x O_R,
\end{equation}
with boundary conditions
%
\begin{equation}
\label{eqHJBRboundarycond} \lambda^R(t,x) = \lambda_1(x)
\qquad\mbox{for all }(t,x) \in\bigl( [0,1) \x\partial O_R \bigr) \cup
\bigl(\{1\} \x O_R\bigr),
\end{equation}
where $\partial O_R$ denotes the boundary of $O_R$.
\end{Theorem}
\begin{pf} First, it follows by the same arguments as in Theorem \ref
{theouviscositysol} that $\lambda^R$ is a viscosity solution of
\eqref{eqHJBR} with boundary condition \eqref{eqHJBRboundarycond}.
The uniqueness follows by the comparison principle of \eqref{eqHJBR},
\eqref{eqHJBRboundarycond}, which holds clearly true from
discussions in Example 3.6 of Crandall et al. \cite{Crandall1992}.
\end{pf}

\subsubsection{Approximation of the transportation cost value}

In the bounded characteristics case, we can give a first approximation
of the minimum transportation cost. Nevertheless, a complete resolution
needs a numerical approximation which will be provided in Section \ref
{secnumericalapprox}. Let us fix the two probability measures $\mu_0$
and $\mu_1$, and simplify the notation $V(\mu_0, \mu_1)$ [resp., $\Vc
(\mu_0, \mu_1)$] to $V$ (resp., $\Vc$).

First, under Assumptions~\ref{assumLconvex}, \ref
{assumLunifcontinuity},~\ref{assumLcroissantp}, \ref
{assumLunifcontinuityinX} and~\ref{assumLbounded}, it follows by
our duality result of Theorem~\ref{theodual} together with Theorem
\ref{theoweakdual} that
\begin{eqnarray}
\label{eqVcVcbequiv} V &=& \Vc:=\sup_{\lambda_1 \in C_b(\R^d)} \bigl(\mu_0(
\lambda_0) - \mu_1(\lambda_1) \bigr)
\nonumber
\\[-8pt]
\\[-8pt]
\nonumber
&=& \overline{\Vc}:=\sup_{\lambda_1 \in C^\infty_b(\R^d)} \bigl(\mu_0(
\lambda_0) - \mu_1(\lambda_1) \bigr),
\nonumber
\end{eqnarray}
where $\lambda_0$ is defined in (\ref{eqlambda0L}).

Let $\operatorname{Lip}^0_K$ denote the collection of all bounded
$K$-Lipschitz-continuous functions $\phi\dvtx \R^d\longrightarrow\R$ with
$\phi(0)=0$, and denote $\operatorname{Lip}^0:=\bigcup_{K>0}\operatorname{Lip}^0_K$. Since
$v(\lambda_1+c)=v(\lambda_1)$ for any $\lambda_1\in C_b(\R^d)$ and
$c\in
\R$, we deduce from \eqref{eqVcVcbequiv} that
\begin{eqnarray*}
V = \sup_{\lambda_1 \in\operatorname{Lip}^0} v(\lambda_1)\qquad \mbox{where } v(
\lambda_1):=\mu_0(\lambda_0) -
\mu_1(\lambda_1).
\end{eqnarray*}
As a first approximation, we introduce the function
%
\begin{equation}
\label{VK} V^K:= \sup_{\lambda_1\in\operatorname{Lip}^0_K} v(\lambda_1).
\end{equation}
Under Assumptions~\ref{assumLbounded} and \ref
{assummu01integrability}, it is clear that $V^K<\infty, \forall K
> 0$ by Lem\-ma~\ref{lemmlambdaLip}. Then, it is immediate that
%
\begin{eqnarray}
\label{VKtoV} \bigl(V^K\bigr)_{K>0} \mbox{ is increasing
and } V^K\longrightarrow V \mbox{ as } K\to\infty.
\end{eqnarray}

Letting $\lambda^R$ be defined in \eqref{eqlambdaR} for every $R >
0$, denote
%
\begin{eqnarray}
\label{eqVKR}  V^{K,R} &:= &\sup_{\lambda_1 \in\operatorname{Lip}^K_0} v^R(
\lambda_1 )\qquad\mbox{and}
\nonumber
\\[-8pt]
\\[-8pt]
\nonumber
 v^R(\lambda_1) &:=& \mu_0\bigl(
\lambda^R_0 \1_{O_R} \bigr)-\mu_1(
\lambda_1 \1_{O_R}).
\end{eqnarray}
Then the second approximation is on variable $R$.
%
\begin{Proposition}
Let Assumptions~\ref{assumLbounded} and \ref
{assummu01integrability} hold true, then for all \mbox{$K>0$},
%
\begin{eqnarray}
\label{eqerrorUKRUK} && \bigl| V^{K,R}-V^K \bigr|
\nonumber
\\[-8pt]
\\[-8pt]
\nonumber
&&\qquad\le C(1+K) \biggl( e^{-R^2/8M + R/2} +\int_{O_{R/2}^c}\bigl(1+|x|\bigr) (
\mu_0 + \mu_1) (dx) \biggr).
\end{eqnarray}
\end{Proposition}
\begin{pf} By their definitions in \eqref{VK} and \eqref{eqVKR}, we have
\begin{eqnarray*}
&&\bigl |V^{K,R}-V^K\bigr |
\\
&&\qquad= \Bigl|\sup_{\lambda_1\in\operatorname{Lip}_K^0} \bigl\{\mu_0\bigl(\lambda_0^R
\1_{O_R}\bigr)-\mu_1(\lambda_1
\1_{O_R}) \bigr\} -\sup_{\lambda_1\in\operatorname{Lip}_K^0} \bigl\{\mu_0(
\lambda_0)-\mu_1(\lambda_1) \bigr\} \Bigr|
\\
&&\qquad\le \sup_{\lambda_1\in\operatorname{Lip}_K^0}\bigl |\mu_0\bigl(\lambda_0^R
\1_{O_R}\bigr)-\mu_0(\lambda_0)\bigr | + K \int
_{O_R^c} |x| \mu_1(dx).
\end{eqnarray*}
Now for all $\lambda_1\in\operatorname{Lip}_K^0$, we estimate from Lemmas
\ref
{lemmlambdaLip} and~\ref{lemmerrorlambdalambdaR} that
\begin{eqnarray*}
&& \bigl|\mu_0\bigl(\lambda_0^R
\1_{O_R}\bigr)-\mu_0(\lambda_0) \bigr|
\\
&&\qquad\le \mu_0 \bigl( \bigl|\lambda_0^R-
\lambda_0 \bigr|\1_{O_{{R}/{2}}} \bigr) + \mu_0 \bigl( \bigl(
\bigl|\lambda_0^R \bigr| + |\lambda_0 | \bigr)
\1_{( O_{{R}/{2}})^c} \bigr)
\\
&&\qquad\le C(1+K) \biggl( \int_{O_{{R}/{2}}}e^{-(R^M_{|x|})^2/2}
\mu_0(dx) + \int_{(O_{{R}/{2}})^c}\bigl(1 + |x| \bigr)
\mu_0(dx) \biggr).
\end{eqnarray*}
Observing that $( R^M_{|x|} )^2 \ge R^2/4M - R + M$ on $O_{{R}/{2}}$, this implies that
\begin{eqnarray*}
&& \bigl|\mu_0\bigl(\lambda_0^R
\1_{O_R}\bigr)-\mu_0(\lambda_0) \bigr|
\\
&&\qquad\le C(1+K) \biggl( e^{-R^2/8M + R/2} + \int_{(O_{{R}/{2}})^c}\bigl(1+|x|\bigr)
\mu_0(dx) \biggr),
\end{eqnarray*}
and the required estimate follows.
\end{pf}

\section{Numerical approximation}
\label{secnumericalapprox}

Throughout this section, we consider the Markovian context where
$L(t,\xb,u) = \ell(t, \xb(t), u)$ under bounded characteristics. Our
objective is to provide an implementable numerical algorithm to compute
$V^{K,R}$ in \eqref{eqVKR}, which is itself an approximation of the
minimum transportation cost $V$ in \eqref{eqVcVcbequiv}.

Although there are many numerical methods for nonlinear PDEs, our
problem concerns the maximization over the solutions of a class of
nonlinear PDEs. To the best of our knowledge, it is not addressed in
the previous literature. In Bonnans and Tan \cite{BonnansTan}, a
similar but more specific problem is considered. Their set $U$ allowing
for unbounded diffusions is out of the scope of this paper. However, by
using the specific structure of their problem, their key observation is
to convert the unconstrained control problem into an optimal stopping
problem for which they propose a numerical approximation scheme. Our
numerical approximation is slightly different, as we avoid the issue of
singular stochastic control by restricting to bounded controls, but
uses their gradient algorithm for the minimization over the choice of
Lagrange multipliers $\lambda_1$.

In the following, we shall first give an overview of the numerical
methods for nonlinear PDEs in Section~\ref{subsecnumerPDE}.
Then by constructing the finite difference scheme for nonlinear PDE
\eqref{eqHJBR}, we get a discrete optimization problem in Section
\ref
{subsecFDapproximation} which is an approximation of $V^{K,R}$.
We then provide a gradient algorithm for the resolution of the discrete
optimization problem in Section~\ref{subsecGradientProjAlgo}.
Finally, we implement our numerical algorithm to test its efficiency in
Section~\ref{subsecnumericalexample}.

In the remaining part of this paper, we restrict the discussion to the
one-dimensional case
\[
d=1\qquad \mbox{so that } O_R = (-R, R).
\]

\subsection{Overview of numerical methods for nonlinear PDEs}
\label{subsecnumerPDE}

There are several numerical schemes for nonlinear PDEs of the form
\eqref{eqHJB}, for example, the finite difference scheme,
semi-Lagrangian scheme and Monte-Carlo schemes. General convergence is
usually deduced by the monotone convergence technique of Barles and
Souganidis \cite{BarlesSouganidis} or the controlled Markov-chain
method of Kushner and Dupuis~\cite{Kushner1990}. Both methods demand
the monotonicity of the scheme, which implies that in practice we
should assume the boundedness of drift and diffusion functions [see,
e.g., the CFL condition \eqref{eqCFLcondition} below]. To derive a
convergence rate, we usually apply Krylov's perturbation method; see,
for example, Barles and Jakobsen \cite{Barles2007}.

For the finite difference scheme, the monotonicity is guaranteed by
the CFL condition [see, e.g., \eqref{eqCFLcondition} below] in the
one-dimensional case $d=1$. However, in the general $d$-dimensional
case, it is usually hard to construct a monotone scheme. Kushner and
Dupuis \cite{Kushner1990} suggested a construction when all covariance
matrices are diagonal dominated. Bonnans et al. \cite
{Bonnans2004,Bonnans2003} investigated this issue and provided an
implementable but
sophisticated algorithm in the two-dimensional case. Debrabant and
Jakobsen \cite{DebrabantJakobsen} proposed recently a semi-Lagrangian
scheme for nonlinear equations of the form \eqref{eqHJB}. However, to
be implemented, it still needs to discretize the space and then to use
an interpolation technique. Therefore, it can be viewed as a kind of
finite difference scheme.

In the high-dimensional case, it is generally preferred to use
Monte-Carlo schemes. For linear and semilinear parabolic PDEs, the
Monte-Carlo methods are usually induced by the Feynman--Kac formula and
backward stochastic differential equations (BSDEs). This scheme is then
generalized by Fahim, Touzi and Warin~\cite{FahimTouziWarin} for
fully nonlinear PDEs. The idea is to approximate the derivatives of the
value function arising in the PDE by conditional expectations, which
can then be estimated by simulation-regression methods. However, the
Monte-Carlo method is not convenient to be used here since for every
terminal condition $\lambda_1$, one needs to simulate many paths of a
stochastic differential equation and then to solve the PDE by
regression method, which makes the computation too costly.

For our problem in \eqref{eqVKR}, we finally choose to use the finite
difference scheme for the resolution of $\lambda^R_0$ since it is easy
to be constructed explicitly as a monotone scheme under explicit
conditions in our context.

\subsection{A finite differences approximation}
\label{subsecFDapproximation}

Let $(l,r) \in\N^2$ and $ h = (\Delta t, \Delta x) \in(\R^+)^2$ be
such that $l \Delta t = 1$ and $r \Delta x = R$. Denote $x_i:= i
\Delta x$, $t_k:= k \Delta t$ and define the discrete grids:
\begin{eqnarray*}
\Nc&:=& \{ x_i \dvtx i \in\Z \}, \qquad \Nc_R:= \Nc\cap(-R,R),
\\
\Mc_{T,R} &:=& \bigl\{ (t_k, x_i) \dvtx(k,i)
\in\Z^+ \x\Z \bigr\} \cap\bigl( [0,1] \x(-R,R) \bigr).
\end{eqnarray*}
The terminal set and boundary set as well as the interior set of $\Mc_{T,R}$ are denoted by
\begin{eqnarray*}
\partial_T \Mc_{T,R}&:=& \bigl\{ (1, x_i)
\dvtx x_i \in\Nc_R \bigr\},\qquad \partial_R
\Mc_{T,R}:= \bigl\{ (t_k, \pm R) \dvtx k=0,\ldots, l
\bigr\},
\\
\mathring{\Mc}_{T,R} &:=& \mathcal{M}_{T,R} \setminus(
\partial_T \mathcal{M}_{T,R} \cup\partial_R
\mathcal{M}_{T,R} ).
\end{eqnarray*}

We shall use the finite differences method to solve the dynamic
programming equation \eqref{eqHJBR}, \eqref{eqHJBRboundarycond} on
the grid $\Mc_{T,R}$. For a function $w$ defined on $\Mc_{T,R}$, we
introduce the discrete derivatives of $w$:
\[
D^{\pm} w(t_k, x_i):= \frac{ w(t_k, x_{i\pm1}) - w(t_k,
x_i)}{\Delta x}
\]
and
\[
(bD)w:=b^+D^+w+b^-D^-w \qquad\mbox{for } b\in\R,
\]
where $b^+:= \max(0, b)$, $b^-:= \max(0, -b)$; and
\[
D^2 w(t_k, x_i):= \frac{ w(t_k, x_{i + 1}) - 2 w(t_k, x_i) + w(t_k, x_{i - 1})} {
\Delta x^2}.
\]
We now define the function $\hat\lambda^{h,R}$(or $\hat\lambda^{h,R,\hat
\lambda_1}$ to emphasize its dependence on the boundary condition
$\hat
\lambda_1$) on the grid $\Mc_{T,R}$ by the following explicit finite
differences approximation of the dynamic programming equation \eqref{eqHJBR}:
\[
\hat\lambda^{h,R}(t_k, x_i) = \hat
\lambda_1 (x_i) \qquad\mbox{on } \partial_T
\Mc_{T,R} \cup\partial_R \Mc_{T,R},
\]
and on $\mathring{\Mc}_{T,R}$,
%
\begin{eqnarray}
\label{eqdiscretwR} && \hat\lambda^{h,R}(t_k,
x_i)
\nonumber
\\[-6pt]
\\[-10pt]
\nonumber
&& \qquad=\biggl(\hat\lambda^{h,R} +\Delta t \inf_{u=(a,b)\in U} \biggl\{
\ell(\cdot, u) + (bD)\hat\lambda^{h,R} + \frac12 a D^2 \hat
\lambda^{h,R} \biggr\} \biggr) (t_{k+1},x_i).\hspace*{-25pt}
\end{eqnarray}
We then introduce the following natural approximation of $v^R$:
%
\begin{equation}
\label{eqvRh}\qquad \hat v^R_h( \hat\lambda_1):= \mu_0 \bigl(\lin^R\bigl[\hat\lambda^{h,R}_0
\bigr] \bigr) -\mu_1 \bigl(\lin^R[\hat
\lambda_1] \bigr)\qquad \mbox{where } \hat\lambda^{h,R}_0:=
\hat\lambda^{h,R}(0,\cdot),
\end{equation}
and for all functions $\phi$ defined on the grid $\Nc_R$ we denote by
$\lin^R[\phi]$ the linear interpolation of $\phi$ extended by zero
outside $[-R,R]$.

We shall also assume that the discretization parameters $h = (\Delta t,
\Delta x)$ satisfy the CFL condition
%
\begin{eqnarray}
\label{eqCFLcondition} \Delta t \biggl( \frac{|b|}{\Delta x} + \frac{|a|}{\Delta x^2}
\biggr) \le1 \qquad \mbox{for all } (a,b)\in U.
\end{eqnarray}
Then the scheme \eqref{eqdiscretwR} is $L^{\infty}$-monotone, so
that the convergence of the scheme is guaranteed by the monotonic
scheme method of Barles and Souganidis \cite{BarlesSouganidis}.
For our next result, we assume that the following error estimate holds.

\begin{Assumption}\label{assumviscosityerrorrate}
There are positive constants $L_{K,R}$, $\rho_1$, $\rho_2$ which are
independent of $h = (\Delta t, \Delta x)$, such that
\begin{eqnarray*}
\mu_0 \bigl( \bigl| \lin^R\bigl[\hat\lambda^{h,R}_0
\bigr]-\lambda_0 \1_{[-R,R]}\bigr | \bigr) \le L_{K,R}
\bigl(\Delta t^{\rho_1} + \Delta x^{\rho_2}\bigr)
\end{eqnarray*}
for all $\lambda_1 \in\operatorname{Lip}^K_0 \mbox{ and } \hat\lambda_1=\lambda_1|_{\Nc_R} $.

\end{Assumption}

Let $\operatorname{Lip}_0^{K,R}$ be the collection of all functions on the grid
$\Nc_R$ defined as restrictions of functions in $\operatorname{Lip}_0^{K}$:
%
\begin{eqnarray}
\label{eqLip0KR} \operatorname{Lip}_0^{K,R} &:=& \bigl\{\hat
\lambda_1:=\lambda_1|_{\Nc_R} \mbox{ for some }
\lambda_1\in\operatorname{Lip}_0^{K} \bigr\}.
\end{eqnarray}
The above approximation of the dynamic value function $\lambda$
suggests the following natural approximation of the minimal
transportation cost value:
%
\begin{eqnarray}
\label{eqVKRh} V^{K,R}_h &:=& \sup_{\hat\lambda_1\in\operatorname{Lip}_0^{K,R}} \hat
v^R_h( \hat\lambda_1)
\nonumber
\\[-8pt]
\\[-8pt]
\nonumber
&=& \sup_{\hat\lambda_1\in\operatorname{Lip}_0^{K,R}} \mu_0 \bigl(\lin^R\bigl[\hat
\lambda^{h,R}_0\bigr] \bigr) -\mu_1 \bigl(
\lin^R[\hat\lambda_1] \bigr).
\end{eqnarray}

\begin{Remark}
Under Assumption~\ref{assumLbounded} and the additional condition
that $\ell$ is uniformly $\frac{1}{2}$-H\"older in $t$ with constant
$M$, then in spirit of the analysis in Barles and Jakobsen \cite
{Barles2007}, Assumption~\ref{assumviscosityerrorrate} holds true
with $\rho_1 = \frac{1}{10}$, $\rho_2 = \frac{1}{5}$ and $L_{K,R} = C(1
+ K +K R)$ with some constant $C$ depending on $M$. This rate is not
the best, but to the best of our knowledge, it is the best rate which
has been proved.
\end{Remark}

\begin{Theorem}
Let Assumption~\ref{assumviscosityerrorrate} be true, then with the
constants $L_{K,R}$, $\rho_1$, $\rho_2$ introduced in Assumption \ref
{assumviscosityerrorrate}, we have
\begin{eqnarray*}
\bigl| V^{K,R}_h - V^{K,R} \bigr| \le L_{K,R}
\bigl(\Delta t^{\rho_1} + \Delta x^{\rho_2}\bigr) + K \Delta x.
\end{eqnarray*}
\end{Theorem}

\begin{pf} First, given $\lambda_1 \in\operatorname{Lip}^K_0$, we take $\hat
\lambda_1:= \lambda_1 |_{\Nc_R} \in\operatorname{Lip}^{K,R}_0$, then clearly $|
\lin^R[\hat\lambda_1] - \lambda_1 |_{L^{\infty}([-R,R])} \le K \Delta x$,
and it follows from Assumption~\ref{assumviscosityerrorrate} and
\eqref{eqVKR} as well as \eqref{eqvRh} that $v^R(\lambda_1) \le
\hat
v^R_h(\hat\lambda_1) + L_{K,R} (\Delta t^{\rho_1} + \Delta x^{\rho_2})
+ K \Delta x $. Hence,
\[
V^{K,R} \le V^{K,R}_h + L_{K,R} \bigl(
\Delta t^{\rho_1} + \Delta x^{\rho_2}\bigr) + K \Delta x.
\]
Next, given $\hat\lambda\in\operatorname{Lip}^{K,R}_0$, let $\lambda_1:=\lin
[\hat\lambda_1] \in\operatorname{Lip}^K_0$ be the linear interpolation of
$\hat
\lambda_1$. It follows from Assumption~\ref{assumviscosityerrorrate}
that $\hat v^R_h(\hat\lambda_1) \le v^R(\lambda_1) + L_{K,R} (\Delta
t^{\rho_1} + \Delta x^{\rho_2})$ and, therefore,
\[
V^{K,R}_h \le V^{K,R} + L_{K,R} \bigl(
\Delta t^{\rho_1} + \Delta x^{\rho_2}\bigr).
\]
\upqed\end{pf}

\subsection{Gradient projection algorithm}
\label{subsecGradientProjAlgo}

In this section we suggest a numerical scheme to approximate $V^{K,R}_h
= \sup_{\hat\lambda_1\in\operatorname{Lip}_0^{K,R}} \hat v^R_h(\hat
\lambda_1) $
in \eqref{eqVKRh}. The crucial observation for our methodology is the
following.
By $B(\Nc_R)$, we denote the set of all bounded functions on $\Nc_R$.

\begin{Proposition} \label{propvconcave}
Under the CFL condition \eqref{eqCFLcondition}, the function $ \hat
\lambda_1 \mapsto\hat v^R_h( \hat\lambda_1)$ is concave on $B(\Nc_R)$.
\end{Proposition}

\begin{pf} Letting $\ub= (\ub_{k,i})_{0 \le k < l, -r < i < r}$, with
$\ub_{k,i} = (\abb_{k,i}, \bbb_{k,i}) \in U$, we introduce $\overline
{\lambda}{}^{h,\ub, \hat\lambda_1}$ (or just $\overline{\lambda
}{}^{h,\ub}$
if there is no risk of ambiguity) as the unique solution of the
discrete linear system on $\Mc_{T,R}$ with a given $\hat\lambda_1$:
\[
\overline{\lambda}{}^{h,\ub}(t_k, x_i) = \hat
\lambda_1(x_i) \qquad\mbox{for } (t_k,x_i)
\in\partial_T \Mc_{T,R} \cup\partial_R
\Mc_{T,R},
\]
and on $\mathring{\Mc}_{T,R}$
%
\begin{eqnarray}
\label{eqdiscretwReb} &&\overline{\lambda}{}^{h,\ub}(t_k,
x_i)
\nonumber
\\[-8pt]
\\[-8pt]
\nonumber
&& \qquad = \bigl(\overline{\lambda}{}^{h,\ub} +\Delta t \bigl( \ell(\cdot,
\ub_{k,i}) + (\bbb_{k,i} D) \overline{\lambda}{}^{h,\ub}
+ \abb_{k,i} D^2 \overline{\lambda}{}^{h,\ub} \bigr)
\bigr) (t_{k+1},x_i).
\nonumber
\end{eqnarray}
Let $\overline{\lambda}{}^{h,\ub}_0:=\overline{\lambda}{}^{h,\ub
}(0,\cdot
)$, and define
\[
\overline{v}{}^{R,\ub}_h(\hat\lambda_1):=
\mu_0 \bigl(\lin^R\bigl[\overline{\lambda}{}^{h,\ub}_0
\bigr] \bigr) -\mu_1 \bigl(\lin^R[\hat
\lambda_1] \bigr).
\]
We claim that
%
\begin{equation}
\label{equhR} \hat{v}_h^R(\hat\lambda_1)
= \inf_{\ub\in U^{l (2r-1)}} \overline{v}_h{}^{R,\ub}(\hat
\lambda_1).
\end{equation}
Indeed, under the CFL condition \eqref{eqCFLcondition}, the finite
difference scheme \eqref{eqdiscretwR} as well as~\eqref
{eqdiscretwReb} are both $L^{\infty}$-monotone in the sense of
Barles and Souganidis~\cite{BarlesSouganidis}. Moreover, the linear
interpolation $\hat\lambda_0 \mapsto\lin^R[\hat\lambda_0]$ is also
monotone. Then taking infimum step by step in \eqref{eqdiscretwR}
and \eqref{eqVKRh} is equivalent to taking infimum globally in~\eqref
{equhR}.

Finally, the concavity of $ \hat\lambda_1 \mapsto\hat v^R_h( \hat
\lambda_1)$ follows from its representation as the infimum of linear
maps in \eqref{equhR}.
\end{pf}

By the previous proposition, $V^{K,R}_h$ consists in the maximization
of a concave function, and a natural scheme to approximate it is the
gradient projection algorithm.

\begin{Remark}\label{remFDexistoptimalcontrol}
Since $U$ is compact by Assumption~\ref{assumLbounded}, then for
every function~$\hat\lambda_1$, we have the optimal control $\hat
{u}(\hat\lambda_1) =  ( \hat{u}_{k,i}(\hat\lambda_1)  )_{0
\le k
< l, -r < i < r}$ such that
%
\begin{eqnarray}
\hat\lambda^{h,R}_0 = \overline{\lambda}{}^{h,\hat u(\hat\lambda_1)}_0
\quad\mbox{and}\quad \hat v_h^R( \hat\lambda_1 ) =
\overline{v}_h{}^{R,\hat u(\hat\lambda_1)}(\hat\lambda_1).
\end{eqnarray}
\end{Remark}

Now we are ready to give the gradient projection algorithm for
$V^{K,R}_h$ in \eqref{eqVKRh}.
Given a function $\varphi\in B(\Nc_R)$, we denote by $P_{\operatorname
{Lip}_0^{K,R}}  ( \varphi )$ the projection of $\varphi$ on
$\operatorname{Lip}_0^{K,R}$, where $\operatorname{Lip}_0^{K,R} \subset B(\Nc_R)$ is
defined in \eqref{eqLip0KR}. Of course, the projection depends on the
choice of the norm equipping
$B(\Nc)$ which in turn has serious consequences on the numerics. We
shall discuss this important issue later.

Letting $\gamma:= (\gamma_n)_{n \ge0}$ be a sequence of positive
constants, we propose the following algorithm:

\begin{Algorithm}
To solve problem \eqref{eqVKRh}:
\begin{itemize}
\item(1) Let $\hat\lambda_1^0:= 0$.
\item(2) Given $\hat\lambda_1^n$, compute the super-gradient $\nabla
\hat v_h^R(\hat\lambda_1^n)$ of $\hat\lambda_1 \mapsto\hat v_h^R(
\hat\lambda_1 )$ at $\hat\lambda_1^n$.
\item(3) Let $\hat\lambda_1^{n+1} = P_{\operatorname{Lip}_0^{K,R}}  (
\hat
\lambda_1^n + \gamma_n \nabla\hat v_h^R(\hat\lambda_1^n)  ) $.
\item(4) Go back to step 2.
\end{itemize}
\end{Algorithm}

In the following, we shall discuss the computation of
super-gradient\break
$\nabla\hat v_h^R(\hat\lambda_1)$, the projection $P_{\operatorname
{Lip}_0^{K,R}}$ as well as the convergence of the above gradient
projection algorithm.

\subsubsection{Super-gradient}

Let $\hat\lambda_1 \in B(\Nc_R)$ be fixed. Then, by Remark \ref
{remFDexistoptimalcontrol}, we may find an optimal control $\hat
{u}(\hat\lambda_1) =  ( \hat{u}_{k,i}(\hat\lambda_1)  )_{0
\le k
< l, -r \le i \le r}$, where\break $\hat{u}_{k,i}(\hat\lambda_1) = ( \hat
{a}_{k,i}( \hat\lambda_1), \hat{b}_{k,i}(\hat\lambda_1)) \in U$, for
system \eqref{equhR}. We then denote by $g^j$ the unique solution of
the following linear system on $\Mc_{T,R}$, for every $-r \le j \le r$:
\[
g^j(t_k, x_i) = \delta_{i,j},
\qquad\mbox{for } (t_k,x_i) \in\partial_T
\Mc_{T,R} \cup\partial_R \Mc_{T,R},
\]
and on $\mathring{\Mc}_{T,R}$,
%
\begin{eqnarray}
\label{eqdiscretsupergradient} g^j(t_k,
x_i) = \bigl(g^j +\Delta t \bigl( \bigl(
\hat{b}_{k,i}( \hat\lambda_1) D \bigr) g^j +
\hat{a}_{k,i}( \hat\lambda_1) D^2
g^j \bigr) \bigr) (t_{k+1},x_i).
\end{eqnarray}
Denote $g^j_0:= g^j(0, \cdot)$ and $\delta_j$ a function on $\Nc_R$
defined by $\delta_j(x_i):= \delta_{i,j}$.

\begin{Proposition}
Let CFL condition \eqref{eqCFLcondition} hold true, then the vector
%
\begin{equation}
\label{eqnablauhR} \nabla\hat v_h^R(\hat
\lambda_1):= \bigl( \mu_0 \bigl( \lin^R
\bigl[g^j_0\bigr] \bigr) - \mu_1 \bigl(
\lin^R[\delta_j] \bigr) \bigr)_{ -r \le j \le r }
\end{equation}
is a super-gradient of $\varphi\in B(\Nc_R) \mapsto\hat
v_h^R(\varphi
) \in\R$ at $\hat\lambda_1$.
\end{Proposition}

\begin{pf} Consider the system \eqref{eqdiscretwReb} introduced in the
proof of Proposition~\ref{propvconcave}. Under the CFL condition
\eqref{eqCFLcondition}, by \eqref{equhR}, we have for every
perturbation $\Delta\hat\lambda_1 \in B(\Nc_R)$,
\[
\hat v_h^R (\hat\lambda_1 + \Delta\hat
\lambda_1) = \overline{v}_h{}^{R,\hat{u}(\hat\lambda_1 + \Delta\hat\lambda_1)} (\hat
\lambda_1 + \Delta\hat\lambda_1) \le \overline{v}_h{}^{R,\hat{u}(\hat\lambda_1)}
(\hat\lambda_1 + \Delta\hat\lambda_1),
\]
which implies that
\begin{eqnarray*}
\hat v_h^R (\hat\lambda_1 + \Delta\hat
\lambda_1) - \hat v_h^R (\hat
\lambda_1) \le \overline{v}_h{}^{R,\hat{u}(\hat\lambda_1)} (\hat
\lambda_1 + \Delta\hat\lambda_1) -\overline{v}_h{}^{R,\hat{u}(\hat\lambda_1)}
(\hat\lambda_1 ).
\end{eqnarray*}
We next observe that for fixed $\hat\lambda_1$, the function $
\varphi
\longmapsto\overline{v}_h{}^{R, \hat{u}(\hat\lambda_1)} (\varphi)$ is
linear, and it follows that
%
\begin{eqnarray}
\label{eqnablauhRexpression} \bigl( \overline{v}_h{}^{R,\hat{u}(\hat\lambda_1)} (
\hat\lambda_1 + \delta_j) -\overline{v}_h{}^{R,\hat{u}(\hat\lambda_1)}
(\hat\lambda_1 ) \bigr)_{ -r \le j \le r }
\end{eqnarray}
is a super-gradient of $\varphi\mapsto\hat v_h^R(\varphi) $ at $\hat
\lambda_1$.
Finally, by \eqref{eqdiscretwReb} and \eqref{eqdiscretsupergradient},
$ g^j(t_k, x_i) = \overline{\lambda}{}^{\hat{u}(\hat\lambda_1), \hat
\lambda_1 + \delta_j}(t_k, x_i) - \overline{\lambda}{}^{\hat{u}(\hat
\lambda_1), \hat\lambda_1}(t_k, x_i),$
where
$\overline{\lambda}{}^{\hat{u}(\hat\lambda_1), \hat\lambda_1 +
\delta_j}$
is the solution of \eqref{eqdiscretwReb} with boundary condition
$\hat\lambda_1 + \delta_j$.
By the definition of
$\overline{v}_h{}^{R, \ub}(\hat\lambda_1)$ in \eqref{equhR}, it
follows that the super-gradient \eqref{eqnablauhRexpression} is
equivalent to
$\nabla\hat v_h^R(\hat\lambda_1)$
defined in~\eqref{eqnablauhR}.
\end{pf}

\subsubsection{Projection}

To compute the projection $P_{\operatorname{Lip}_0^{K,R}}  ( \varphi )$,
$\forall\varphi\in B(\Nc_R)$, we need to equip $B(\Nc_R)$ with a
specific norm. In order to obtain a simple projection algorithm, we
shall introduce an invertible linear map between $B(\Nc_R)$ and $\R^{2r+1}$,
then equip on $B(\Nc_R)$ the norm induced by the classical
$L^2$-norm on $\R^{2r+1}$.

Let us define the invertible linear map $\Tc_R$ from $B(\Nc_R)$ to
$\R^{2r+1}$ as
\begin{eqnarray*}
\psi_i = \Tc_R (\varphi)_i:= \cases{
\varphi(x_{i+1}) - \varphi(x_i),&\quad $i = 1, \ldots, r,$
\vspace*{2pt}
\cr
\varphi(0),&\quad $i = 0,$\vspace*{2pt}
\cr
\varphi(x_{i-1}) -
\varphi(x_i),&\quad $i = -1, \ldots, -r,$ }
\end{eqnarray*}
and define the norm $| \cdot|_R$ on $B(\Nc_R)$ (easily be verified) by
\[
| \varphi|_R:= \bigl| \Tc_R(\varphi) \bigr|_{L^2(\R^{2r+1})}\qquad
\forall\varphi\in B(\Nc_R).
\]
Notice that
\begin{eqnarray*}
\Tc_R \operatorname{Lip}_0^{K,R} &:=& \bigl\{ \psi=
\Tc_R \varphi\dvtx \varphi\in\operatorname{Lip}_0^{K,R}
\bigr\}
\\
&=& \bigl\{ \psi= (\psi_i)_{-r \le i \le r} \in[-K\Delta x, K \Delta
x]^ {2r+1} \dvtx\psi_0 = 0 \bigr\}.
\end{eqnarray*}
Then the projection $P_{\operatorname{Lip}_0^{K,R}}$ from $B(\Nc_R)$ to
$\operatorname
{Lip}_0^{K,R}$ under norm $| \cdot|_R$ is equivalent to the projection
$P_{\Tc_R \operatorname{Lip}_0^{K,R}}$ from $\R^{2r+1}$ to $\Tc_R \operatorname
{Lip}_0^{K,R}$ under the $L^2$-norm, which is simply written as
\[
\bigl(P_{\Tc_R \operatorname{Lip}_0^{K,R}} (\psi ) \bigr)_i = \cases{ 0,& \quad$\mbox{if } i
= 0,$\vspace*{2pt}
\cr
(K \Delta x) \wedge\psi_i \vee(-K \Delta x),&\quad
$\mbox{otherwise}.$}
\]

\subsubsection{Convergence rate}

Now, let us give a convergence rate for the above gradient projection
algorithm. In preparation, we first provide an estimate for the norm of
super-gradients $\nabla\hat v_h^R$.

\begin{Proposition}
Suppose that CFL condition \eqref{eqCFLcondition} holds true, then
$| \hat v_h^R (\varphi_1) - \hat v_h^R (\varphi_2)| \le2 |\varphi_1 -
\varphi_2|_{\infty}$ for every $\varphi_1, \varphi_2 \in B(\Nc_R)$. In
particular, the super-gradient $\nabla\hat v_h^R$ satisfies
%
\begin{equation}
\label{eqSuperGradBound} \bigl| \nabla\hat v_h^R (\hat
\lambda_1) \bigr|_R \le 2 \sqrt{\frac{R}{\Delta x} +1 },
\qquad\mbox{for all } \hat\lambda_1 \in B(\Nc).
\end{equation}
\end{Proposition}
\begin{pf} Under the CFL condition, the scheme \eqref{eqdiscretwR} is
$L^{\infty}$-monotone, then $  | \hat\lambda^{h,R, \varphi_1}_0 -
\hat\lambda^{h,R, \varphi_2}_0  |_{\infty} \le|\varphi_1 -
\varphi_2|_{\infty}, $ and it follows from the definition of $\hat v_h^R$ in~\eqref{eqvRh} that
%
\begin{equation}
\label{eqgraduhRinfty} \bigl| \hat v_h^R (
\varphi_1) - \hat v_h^R (
\varphi_2)\bigr| \le 2 |\varphi_1 - \varphi_2|_{\infty}.
\end{equation}
Next, by the Cauchy--Schwarz inequality,
\begin{eqnarray*}
|\varphi_1 - \varphi_2|_{\infty} &\le& \max
\Biggl( \sum_{i=0}^r \bigl|\Tc_R(
\varphi_1 - \varphi_2)_i\bigr |, \sum
_{i=0}^{-r} \bigl|\Tc_R(\varphi_1
- \varphi_2)_i\bigr | \Biggr)
\\
&\le& \sqrt{r+1} |\varphi_1 - \varphi_2|_R.
\end{eqnarray*}
Together with \eqref{eqgraduhRinfty}, this implies that \eqref
{eqSuperGradBound} holds for every super-gradient $\nabla\hat v_h^R
(\hat\lambda_1)$.
\end{pf}

Let us finish this section by providing a convergence rate for our
gradient projection algorithm. Denote
\[
\Pi:= \max_{\varphi_1, \varphi_2 \in\operatorname{Lip}_0^{K,R} } | \varphi_1 - \varphi_2
|_R^2 \le 2r (K \Delta x)^2 \le2
K^2 R \Delta x,
\]
and it follows from Section 5.3.1 of Ben-Tal and Nemirovski \cite
{BenTal2001} that
%
\begin{eqnarray}
\label{eqConvRateGradientProj} 0 \le V_h^{K,R} -
\max_{n \le N} \hat v_h^R \bigl(\hat
\lambda_1^n\bigr) &\le& \frac{ \Pi+ \sum_{n=1}^N \gamma_n^2 | \nabla\hat v_h^R(\hat
\lambda_1^n) |_R^2 }{\sum_{n=1}^N \gamma_n}
\nonumber
\\[-8pt]
\\[-8pt]
\nonumber
&\le& \frac{ 2 K^2 R \Delta x + 4  ( {R}/{\Delta x} +1  )
\sum_{n=1}^N \gamma_n^2 }{\sum_{n=1}^N \gamma_n}.
\end{eqnarray}

We have several choices for the series $\gamma= (\gamma_n)_{n \ge1}$:
\begin{itemize}
\item Divergent series: $ \gamma_n \ge0$, $\sum_{n=1}^{\infty}
\gamma_n = +\infty$ and $\sum_{n=1}^{\infty} \gamma_n^2 < +\infty$, then the
right-hand side of \eqref{eqConvRateGradientProj} converges to $0$ as
$N \to\infty$.
\item Optimal stepsizes: $\gamma_n = \frac{ \sqrt{2 \Pi} } { |
\nabla
\hat v_h^R (\hat\lambda_1^n) |_R \sqrt{n} } $, \cite{BenTal2001}
shows that
\begin{eqnarray*}
V_h^{K,R} - \max_{n \le N} \hat
v_h^R \bigl(\hat\lambda_1^n
\bigr) &\le& C_1 \frac{ ( \max_{1 \le n \le N} | \nabla\hat v_h^R( \hat\lambda_1^n) |_R ) \cdot\sqrt{2 \Pi}}{\sqrt{N}}
\\
&\le& C \frac{K ( R+\sqrt{R \Delta x})}{\sqrt{N}}
\end{eqnarray*}
for some constant $C$ independent of $K$, $R$, $\Delta t$, $\Delta x$
and $N$.
\end{itemize}

\subsection{Numerical examples}
\label{subsecnumericalexample}

We finally implement the above algorithm in the context of an
application in finance which consists in the determination of the
optimal no-arbitrage bounds of exotic options.

As discussed in the \hyperref[secintroduction]{Introduction}, this
problem has been solved by means of the Skorokhod Embedding Problem
(SEP) in the context of some specific examples of derivative
securities. However, the SEP approach is not suitable for numerical
approximation. In Davis, Obloj and Raval \cite{DavisOblojRaval}, the
authors consider a similar problem for the weighted variance swap
option which can be included in our context. In contrast to our
constraint $\P\circ X_1^{-1} = \mu_1$ in~\eqref{Pmu0mu1}, they impose
the constraint of the form $\E^{\P} [\phi_k(X_1)] = p_k, k=1,\ldots,
n$ for some functions $(\phi_k)_{1 \le k \le n}$ and constants
$(p_k)_{1 \le k\le n}$. Then, they convert their problem into a
semi-infinite linear programming problem which can be solved
numerically. We shall use some techniques in~\cite{DavisOblojRaval}
to derive an explicit solution for some examples in order to compare
with our numerical results.

\subsubsection{A toy example}

\begin{figure}[b]

\includegraphics{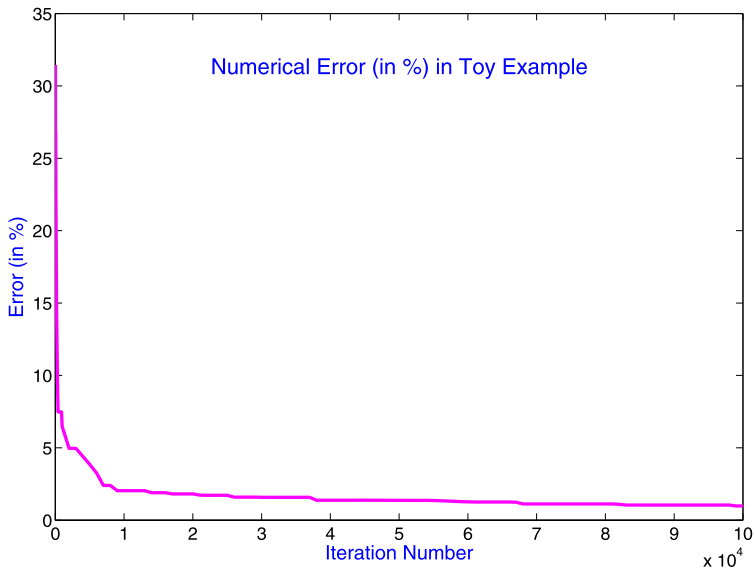}

\caption{Numerical Example 1 (toy example): $\mu_i = N(0, \sigma_i^2)$
with $\sigma_0 = 0.1$, $\sigma_1 = 0.2$, $K = 1.5$, $R = 1$, $\Delta x
= 0.1$, $\Delta t = 0.025$. The computation time is 56.38 seconds for
$10^5$ iterations. }
\label{figToyExample}
\end{figure}

Suppose that $\ell(t,x,a,b) = a$, and $U:= [\underline{a}, \overline
{a}] \x\{ 0 \}$, then under every $\P\in\Pc$, the canonical process
$X$ is a martingale. Suppose that $\Pc(\mu_0, \mu_1)$ is nonempty, then
it is clear that
\begin{eqnarray*}
V &=& \inf_{\P\in\Pc(\mu_0, \mu_1)} \E^{\P} \biggl[ \int_0^1
\alpha^{\P
}_t\, dt \biggr]
\\
&=& \inf_{\P\in\Pc(\mu_0, \mu_1)} \E^{\P} \bigl[ X^2_1
- X^2_0 \bigr] = \mu_1(\phi_0) -
\mu_0(\phi_0),
\end{eqnarray*}
where $\phi_0(x):= x^2$. In our implemented example, we choose $\mu_i$
as normal distribution $N(0, \sigma_i^2)$ with $\sigma_0 = 0.1$,
$\sigma_1 = 0.2$, $\underline{a} = 0$ and $\overline{a} = 0.1$. It follows by
direct computation that $V = 0.03$. In our numerical test, for $10^5$
iterations, the computation time is 56.38 seconds, and it gives a
numerical solution 0.029705, which implies that the relative error is
less than $1 \%$; see Figure~\ref{figToyExample}.

\subsubsection{The weighted variance swap contract}

Let $S = (S_t)_{t \ge0}$ denote the price process of an underlying
stock. We assume that $S$ is a scalar positive continuous
semimartingale. The variance swap contract is therefore defined by the
payoff $ \langle\log S \rangle_1$ at maturity $1$, which is the
quadratic variation of the process $(\log S_t)_{ t\ge0}$ at time $1$.

Following Section 4 of \cite{DavisOblojRaval}, we shall consider an
$\eta$-weighted variance swap, for some Lipschitz function $\eta\dvtx  \R
\to\R$. This is a derivative security defined by the payoff at
maturity $1$:
\[
\int_0^1 \eta(\log S_t) d \langle
\log S \rangle_t.
\]
Under no additional information, any martingale measure $\P$ (i.e., a
probability measure under which the process $S$ is a martingale)
induces an admissible no-arbitrage price
\[
\E^{\P} \biggl[ \int_0^1 \eta\bigl(
\log(S_t)\bigr) d \langle\log S \rangle_t \biggr].
\]
Following Galichon et al. \cite{ght}, we assume that all European
options maturing at time $1$ with all possible strikes are liquids and
available for trading, that is, $c_1(y):=\E[(S_1-y)^+]$ is given for
all $y\ge0$. Then the marginal distribution of $S_1$ under $\P$ is
given by $\tilde{\mu}_1[y,\infty)=-\partial^-c_1(y)$. In other words,
for every $\lambda_1 \in C_b(\R)$, the derivative security with payoff
$\lambda_1(S_1)$ at maturity $1$ is available for trading (long or
short) at the no-arbitrage price $\tilde{\mu}_1(\lambda_1)$. Under this
additional information, a no-arbitrage lower bound of the $\eta
$-weighted variance swap is given by
%
\begin{equation}
\label{eqBoundVS} \sup_{\lambda_1 \in C_b(\R)} \biggl\{ \inf_{\P}
\E^{\P} \biggl[ \int_0^1 \eta(\log
S_t)\, d \langle\log S \rangle_t + \lambda_1(S_1)
\biggr] - \tilde{\mu}_1(\lambda_1) \biggr\},
\end{equation}
where the infimum is taken over all martingale measures for $S$.

\begin{figure}

\includegraphics{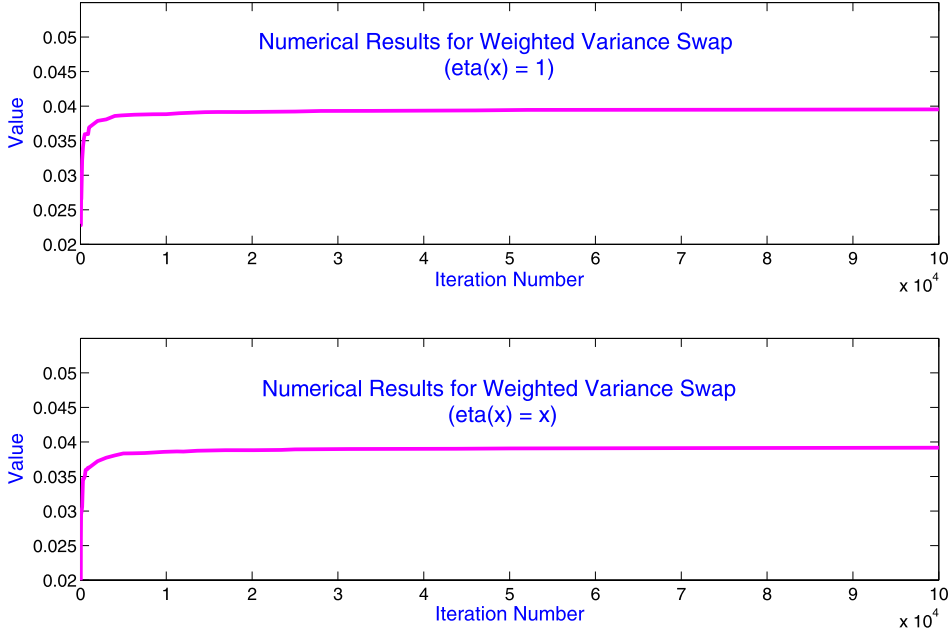}

\caption{Numerical Example 2 (weighted variance swap): $\sigma= 0.2$,
$K = 1.5$, $R = 2$, $\Delta x = 0.1$, $\Delta t = 0.025$. For weight
function $\eta_1(x) = 1$, the numerical solution is 0.0395311 after
$10^5$ iterations. For weight function $\eta_2(x) = x$, the numerical
solution is 0.0391632 after $10^5$ iterations.}
\label{figVarianceSwap}
\end{figure}

This problem can be studied by our mass transportation problem. Suppose
that $S_t = \exp(X_t)$, where $X$ is the canonical process on the
canonical space $\O$. Suppose further that
\[
U:= \bigl\{ \bigl( a, -\tfrac{1}{2}a \bigr) \in S_1 \x\R
\dvtx a \in [\underline{a}, \overline{a}] \bigr\},
\]
with positive constants $\underline{a} \le\overline{a} < \infty$. Then
under every $\P\in\Pc$ [defined below~\eqref
{eqconvexitysetSigma}], the process $S_t:= \exp(X_t)$ is a positive
continuous martingale. If we take the infimum in \eqref{eqBoundVS}
over $\Pc$, it follows by our duality result (Theorem~\ref{theodual})
that the bound \eqref{eqBoundVS} equals
\[
V = \inf_{\P\in\Pc(\delta_{x_0}, \mu_1)} \E^{\P} \int_0^1
\eta(X_t) \a^{\P}_t \,dt,
\]
where $x_0 = \log S_0 \in\R$ and $\mu_1$ is the distribution of $X_1 =
\log S_1$ when $S_1 \sim\tilde\mu_1$ [or, equivalently, it is derived
from $\tilde\mu_1$ by $ \int_{\R} \varphi(x) \mu_1(dx):= \int_{\R}
\varphi(\log y) \tilde\mu_1(dy)$, $\forall\varphi\in C_b(\R)$].
Furthermore, by similar techniques as in Section 4 of Davis et~al.~\cite
{DavisOblojRaval}, using It\^o's formula, it follows that when $\Pc
(\delta_{x_0}, \mu_1)$ is nonempty, we have
%
\begin{eqnarray}
\label{eqBoundVS2} V = \inf_{\P\in\Pc(\delta_{x_0}, \mu_1)} \E^{\P} \bigl[
\phi(X_1) - \phi(X_0)\bigr] = \mu_1(\phi) -
\phi(x_0), \label{varianceswap-explicit}
\end{eqnarray}
where $\phi$ is a solution to $\phi''(x) - \phi'(x) = 2 \eta(x)$.

In our numerical experiments, we choose $\underline{a} = 0$,
$\overline
{a} = 0.1$, $x_0 = 1$, $\mu_1$ is a normal distribution $N(1 - a/2, a)$
with $a = 0.04 \in[0, 0.1]$. Then $\Pc(\delta_{x_0}, \mu_1)$ is
nonempty since the probability $\P$ induced by the process $(1-at/2 +
\sqrt{a} W_t)_{0 \le t \le1}$ (with Brownian motion $W$) belongs to
it. In a first example, we choose $\eta_1(x) = 1$, then $\phi_1(x):=
-2x + C_1 e^x + C_2$ is the solution to $\phi''(x) - \phi'(x) = 2
\eta_1(x)$. It follows by direct computation that the value in \eqref
{eqBoundVS2} is given by $V_1 = 0.04$. Our numerical solution is
$0.0395311$ after $10^5$ iterations, which takes 138.51 seconds. In a
second example, we choose $\eta_2(x) = x$, then $\phi_2(x):= -x^2 - 2x
- C_1 e^x + C_2$ is the solution to $\phi''(x) - \phi'(x) = 2 \eta_2(x)$.
It follows that the value in \eqref{eqBoundVS2} is given by
$V_2 = a - a^2/4 = 0.0396$. In our numerical test, the computation time
is 142.23 seconds for $10^5$ iterations and it gives the numerical
solution $0.0391632$; see Figure~\ref{figVarianceSwap}.

\begin{appendix}\label{app}

\section*{Appendix}

We first give a result which follows directly from the measurable
selection theorem.
Let $E$ and $F$ be two Polish spaces with their Borel $\sigma$-fields
$\Ec:= \Bc(E)$ and $\Fc:= \Bc(F)$. $A \in\Ec\otimes\Fc$ is a
measurable subset in the product space $E \x F$ satisfying that for
every $x \in E$, there is $y \in F$ such that \mbox{$(x,y) \in A$}. Letting
$\mu$ be a probability measure on $(E, \Ec)$, we denote by $\Ec^{\mu}$
the $\mu$-completed $\sigma$-field of $\Ec$. Suppose that $f \dvtx A
\to
\R
\cup\{ \infty\}$ is $\Ec\otimes\Fc$-measurable, and denote
%
\setcounter{equation}{0}
\begin{equation}
\label{eqdgsupf} g(x):= \inf\bigl\{ f(x, y), (x,y) \in A \bigr\}.
\end{equation}

\setcounter{Theorem}{0}
\begin{Theorem} \label{TheoMeasurabilitySup}
The function $g$ is $\Ec^{\mu}$-measurable. Moreover, for every $\eps
> 0$, there is a $\Ec^{\mu}$-measurable variable $Y_{\eps}$ such that
for $\mu$-a.e. $x \in E$, $(x, Y_{\eps}(x)) \in A$ and
%
\begin{equation}
\label{eqepsmeasurablesection} f\bigl(x, Y_{\eps}(x)\bigr) \le \bigl( g
(x) + \eps \bigr) \1_{g(\om) > -
\infty
} - \frac{1}{\eps} \1_{g(x) = -\infty}.
\end{equation}
\end{Theorem}

\begin{Remark}
Theorem~\ref{TheoMeasurabilitySup} is almost the same as Proposition
7.50 of Bertsekas and Shreve \cite{Bertsekas1978}, and can be easily
deduced from it. The key argument is the measurable selection theorem.
We also refer to Section 12.1 of Stroock and Varadhan \cite
{Stroock1979}, Chapter 7 of Bertsekas and Shreve \cite{Bertsekas1978}
or Chapter 3 of Dellacherie and Meyer \cite{Dellacherie1978} for
different versions of the measurable selection theorem.
\end{Remark}

We next report a theorem which provides the unique canonical
decomposition of a continuous semimartingale under different
filtrations. In particular, it follows that an It\^o process has a
diffusion representation, by taking the filtration generated by itself.
This is in fact Theorem 7.17 of Liptser and Shiryayev \cite
{Liptser1977} in the \mbox{1-dimensional} case or Theorem 4.3 of Wong \cite
{Wong1971} in the multidimensional case.

\begin{Theorem}\label{theoIto2Diffusion}
In a filtrated space $(\O, \F= (\Fc_t)_{0 \le t\le1}, \P)$ (here
$\O$
is not necessarily the canonical space), a process $X$ is a continuous
semimartingale with canonical decomposition:
\[
X_t = X_0 + B_t + M_t,
\]
where $B_0 = M_0 =0,$ and $B = (B_t)_{0 \le t \le1}$ is of finite
variation and $M = (M_t)_{0 \le t \le1}$ a local martingale. In
addition, suppose that there are measurable and $\F$-adapted processes
$(\a,\beta)$ such that
\[
B_t = \int_0^t
\beta_s \,ds,\qquad \int_0^1 \E\bigl[|
\beta_s|\bigr] \,ds < \infty \quad\mbox {and}\quad A_t:= \langle M
\rangle_t = \int_0^t
\a_s \,ds.
\]
Let $\F^X = (\Fc_t^X)_{0 \le t \le1}$ be the filtration generated by
process $X$ and $\bar{\F} = (\bar{\Fc}_t)_{0 \le t \le1}$ be a
filtration such that $\Fc^X_t \subseteq\bar{\Fc}_t \subseteq\Fc_t$.
Then $X$ is still a continuous\vadjust{\goodbreak} semimartingale under $\bar{\F}$, whose
canonical decomposition is given by
\[
X_t = X_0 + \int_0^t
\bar{\beta}_s \,ds + \bar{M}_t \qquad\mbox{with } \bar
{A}_t:=\langle\bar{M} \rangle_t = \int
_0^t \bar{\a}_s \,ds,
\]
where
\[
\bar{\beta}_t = \E[\beta_t | \bar{\Fc}_t]
\quad\mbox{and}\quad \bar{\a}_t = \a_t, d\P\x \,dt\mbox{-a.e.}
\]
\end{Theorem}

\begin{pf*}{Proof of Theorem~\ref{theouviscositysol}} The
characterization of the value function as viscosity solution to a
dynamic programming equation is a natural result of the dynamic
programming principle. Here, we give a proof, similar to that of
Corollary~5.1 in \cite{Bouchard2011}, in our context.

(1) We first prove the subsolution property. Suppose that $(t_0, x_0)
\in[0,1) \x\R^d$ and $\phi\in C_c^{\infty}([0,1) \x\R^d)$ is a
smooth function such that
\[
0 = (\lambda-\phi) (t_0, x_0) > (\lambda- \phi) (t,x)\qquad
\forall(t,x) \neq(t_0, x_0).
\]
By adding $\eps (|t-t_0|^2 + |x-x_0|^4  )$ to $\phi(t,x)$, we
can suppose that
%
\begin{equation}
\label{eqvarphiminusu} \phi(t,x) \ge\lambda(t,x) + \eps \bigl(|t-t_0|^2
+ |x-x_0|^4 \bigr)
\end{equation}
without losing generality. Assume to the contrary that
\[
- \partial_t \phi(t_0, x_0) - H
\bigl(t_0,x_0,D_x \phi(t_0,
x_0), D^2_{xx} \phi(t_0,
x_0) \bigr) > 0,
\]
and we shall derive a contradiction. Indeed, by definition of $H$,
there is $c > 0$ and $(a,b) \in U$ such that
\begin{eqnarray}
- \partial_t \phi(t, x) - b \cdot D_x \phi(t, x) -
\tfrac{1}{2} a \cdot D^2_{xx} \phi(t, x) -
\ell(t,x,a,b) > 0\nonumber\\
\eqntext{\forall(t,x) \in B_c(t_0,
x_0),}
\end{eqnarray}
where $B_c(t_0, x_0):= \{ (t,x) \in[0,1) \x\R^d \dvtx|(t,x) - (t_0,
x_0)| \le c \}$. Let $\tau:= \inf\{ t \ge t_0 \dvtx(t,X_t) \notin
B_c(t_0, x_0) \} \land T$, then
\begin{eqnarray*}
\lambda(t_0, x_0) = \phi(t_0,
x_0) &\ge& \inf_{\Pb\in\Pcb
_{t_0,x_0,0,0}} \E^{\Pb} \biggl[ \int
_t^{\tau} \ell(s,X_s,
\nu_s) \,ds + \phi(\tau, X_{\tau}) \biggr]
\\
&\ge& \inf_{\Pb\in\Pcb_{t_0,x_0,0,0}} \E^{\Pb} \biggl[ \int_t^{\tau}
\ell(s,X_s, \nu_s) \,ds + \lambda(\tau,
X_{\tau}) \biggr] + \eta,
\end{eqnarray*}
where $\eta$ is a positive constant by \eqref{eqvarphiminusu} and
the definition of $\tau$. This is a contradiction to Proposition \ref
{propdpp}.

(2) For the supersolution property, we assume to the contrary that
there is $(t_0, x_0) \in[0,1) \x\R^d$ and a smooth function $\phi$
satisfying
\[
0 = (\lambda-\phi) (t_0, x_0) < (\lambda- \phi) (t,x)\qquad
\forall(t,x) \neq(t_0, x_0)
\]
and
\[
- \partial_t \phi(t_0, x_0) - H \bigl(
t_0,x_0,D_x \phi(t_0,
x_0), D^2_{xx} \phi(t_0,
x_0) \bigr) < 0.
\]
We also suppose without losing generality that
%
\begin{eqnarray}
\label{eqvarphiplusu} \phi(t,x) \le\lambda(t,x) - \eps \bigl(|t-t_0|^2
+ |x-x_0|^4 \bigr).
\end{eqnarray}

By continuity of $H$, there is $c>0$ such that for all $(t,x) \in
B_c(t_0, x_0)$ and every $(a,b) \in U$,
\begin{eqnarray*}
- \partial_t \phi(t, x) - b \cdot D_x \phi(t, x) -
\tfrac{1}{2} a \cdot D^2_{xx} \phi(t, x) -
\ell(t,x,a,b) < 0.
\end{eqnarray*}
Let $\tau:= \inf\{ t \ge t_0 \dvtx(t,X_t) \notin B_c(t_0, x_0) \}
\land T$, then
\begin{eqnarray*}
\lambda(t_0, x_0) = \phi(t_0,
x_0) &\le& \inf_{\Pb\in\Pcb
_{t_0,x_0,0,0}} \E^{\Pb} \biggl[ \phi(\tau,
X_{\tau}) + \int_{t_0}^{\tau
} \ell(s,
X_s, \nu_s) \,ds \biggr]
\\
&\le& \inf_{\Pb\in\Pcb_{t_0,x_0,0,0} } \E^{\Pb} \biggl[ \lambda (\tau,
X_{\tau}) + \int_{t_0}^{\tau} \ell(s,
X_s, \nu_s) \,ds \biggr] - \eta
\end{eqnarray*}
for some $\eta> 0$ by \eqref{eqvarphiplusu}, which is a
contradiction to Proposition~\ref{propdpp}. 
\end{pf*}
\end{appendix}

\section*{Acknowledgments}

We are grateful to J. Fr\'ed\'eric Bonnans, Nicole El Karoui and
Alfred Galichon for fruitful discussions and two anonymous referees for
helpful comments and suggestions.

%

%


\printaddresses

\end{document}